\documentclass[10pt]{amsart}
\usepackage{amssymb}
\usepackage[all]{xy}
\expandafter\edef\csname :RestoreCatcodes\endcsname{%
   \catcode`\noexpand :=\the\catcode`:%
   \catcode`\noexpand @=\the\catcode`@%
   \catcode`\noexpand /=\the\catcode`/%
   \catcode`\noexpand &=\the\catcode`&%
   \catcode`\noexpand \^^M=\the\catcode`\^^M%
   \catcode`\noexpand \^^I=\the\catcode`\^^I%
   \expandafter\let\csname:RestoreCatcodes\endcsname=\noexpand\undefined}
\catcode`\:11  \catcode`\@11
   \let\:wlog\wlog \def\wlog#1{}
   \def\:wrn#1#2{\immediate\write\sixt@@n{--DraTeX warning--
      \ifcase #1
    DraTex.sty already loaded
\or \string\Draw\space within \string\Draw
\or Changing definition of \string#2
\or No intersection points: #2
\or Improper rotation of axes: #2
\or (#2) in \string\DSeg\space is a point
\fi}}
\def\:err#1#2{\errmessage{--DraTeX error-- \ifcase #1
     \string#2\space meaningless in three dimensions
\or  \string#2\space meaningless in two dimensions
\or  No \string\MarkLoc(#2)
\or  \string#2 in three dimensions
\or  Too many parameters in definition
\or  \string\MoveFToOval(#2)?
\fi}}
   \ifx\:Xunits\:undefined \else \:wrn0{} \fi
   \catcode`\ 9  \catcode`\^^M9 \catcode`\^^I9
      \newdimen\:LBorder \newdimen\:RBorder\chardef\:eight=8
\mathchardef\:cccvx=360
\newdimen\:mp    \:mp   0.1\p@
\newdimen\:mmp   \:mmp  0.01\p@
\newdimen\:mmmp  \:mmmp 0.001\p@
\newdimen\:XC    \:XC   90\p@
\newdimen\:CVXXX \:CVXXX180\p@
\newdimen\:CCCVX \:CCCVX\:cccvx\p@ \newdimen\:TeXLoc
\newbox\:box\newif\if:IIID  \newdimen\:Z   \newdimen\:Zunits
\newdimen\:Ex   \newdimen\:Ey  \newdimen\:Ez

\def\:AbsVal#1{ \ifdim#1<\z@-\fi #1 }
\def\:abs#1{\ifdim #1<\z@ #1-#1 \fi}
\def\:AbsDif#1#2#3{  #1#2   \advance#1  -#3
   \ifdim #1<\z@ #1-#1 \fi}
\def\:diff#1#2#3{ #1#2  \advance#1 -#3 }
\def\:average#1#2#3{
   #1#2  \advance#1  #3   \divide#1 \tw@}\def\:Opt#1#2#3#4{
   \def\:temp{
      \ifx      \:next\ifnum \def\:next{#3#1#4#2}
      \else\ifx \:next#1     \def\:next{#3}
      \else                  \def\:next{#3#1#4#2}\fi\fi \:next}
   \futurelet\:next\:temp}\def\Define#1{\:multid#1
   \:Opt(){\:Define#1}0}

\def\:DraCatCodes{\catcode`\ 9   \catcode`\^^M9
   \catcode`\^^I9  \catcode`\&13  \catcode`\~13 }

\def\:Define#1(#2){\begingroup  \:DraCatCodes  \::Define#1(#2)}

\def\::Define#1(#2)#3{\endgroup
   \let\:NextDefine\NextDefine
   \let\NextDefine\relax
   \ifcase#2\relax
      \def#1{#3}\or
      \:TxtPar\def#1(##1){#3}\or
         \:TxtPar\def#1(##1,##2){#3}\or
   \:TxtPar\def#1(##1,##2,##3){#3}\or
   \:TxtPar\def#1(##1,##2,##3,##4){#3}\or
   \:TxtPar\def#1(##1,##2,##3,##4,##5){#3}\or
   \:TxtPar\def#1(##1,##2,##3,##4,##5,##6){#3}\or
   \:TxtPar\def#1(##1,##2,##3,##4,##5,##6,##7){#3}\or
   \:TxtPar\def#1(##1,##2,##3,##4,##5,##6,##7,##8){#3}\or
      \:TxtPar\def#1(##1,##2,##3,##4,##5,##6,##7,##8,##9){#3}\or
      \:err4{}\fi      \let\:TxtPar\relax  \:NextDefine}

\let\NextDefine\relax\let\:TxtPar\relax
\def\WarningOn{\def\:multid##1{
   \ifx ##1\:undefined \else \:wrn2##1\fi}}
\def\:gobble#1{}
\def\WarningOff{\let\:multid\:gobble}     \WarningOff
\Define\Indirect{\futurelet\:next\:Indirect}

\Define\:Indirect{\:theDoReg
   \ifx \:next<
      \def\:temp{\let\DoReg\:DoReg}
      \def\:next<##1>{\expandafter\:temp\csname :<##1>\endcsname}
   \else
      \def\:next##1<##2>{
         \expandafter\ifx \csname :<##2> \endcsname \relax
               \def\:next{##1}     \fi
\:indrwrn\Define     \:indrwrn\Object
\:indrwrn\Table      \:indrwrn\IntVar     \:indrwrn\DecVar
         \def\:temp{\let\DoReg\:DoReg##1}
         \expandafter\:temp \csname :<##2> \endcsname}
   \fi      \:next}  \def\:indrwrn#1{    \def\:temp{#1}
   \ifx \:next\:temp \def\:wrn##1##2{\let\:wrn\::wrn} \fi}
\let\::wrn\:wrn
\Define\:Hline{
   \setbox\:box\hbox{\vrule height0.5\:thickness
                    depth 0.5\:thickness width\:x}
   {\:d\:X \advance\:d \wd\:box
 \advance\:X -\:TeXLoc   \global\:TeXLoc\:d
 \vrule width\:X depth\z@ height\z@
 \raise \:Y \box\:box} }\Define\:Vline{
   \setbox\:box\hbox{\vrule width\:thickness
                 \ifdim \:y>\z@ height\:y  depth\z@
                 \else          height\z@  depth-\:y \fi}
   \advance\:X  -0.5\:thickness
   {\:d\:X \advance\:d \wd\:box
 \advance\:X -\:TeXLoc   \global\:TeXLoc\:d
 \vrule width\:X depth\z@ height\z@
 \raise \:Y \box\:box} }\Define\:MvTo(2){\:X#1\:Xunits \:Y#2\:Yunits}
\Define\:Mv(2){\advance\:X  #1\:Xunits
               \advance\:Y  #2\:Yunits}
\def\:DLn(#1,#2,{\:MvTo(#1,#2) \:LnTo(}\Define\:LnTo(2){
   \:x#1\:Xunits \advance\:x  -\:X
   \:y#2\:Yunits \advance\:y  -\:Y
   \:Ln(\:x\du,\:y\du) }\Define\:Ln(2){
   \:x#1\:Xunits \:y#2\:Yunits
   { \ifdim \:x<\z@
        \advance\:X \:x   \:x-\:x
        \advance\:Y \:y   \:y-\:y
     \fi
     \:yy\:AbsVal\:y
     \:dd\:yy \advance\:dd \:x
     \ifdim \:dd>\:mmmp
        \ifdim \:x>\:yy
           { \ifdim\:X<\:LBorder \global\:LBorder\:X\fi
\advance \:X  \:x
\ifdim \:X>\:RBorder \global\:RBorder\:X\fi }
           \let\:Yunitsy\:x  \let\:Xunitsx\:y  \let\:temp\:Hline
           \:dd0.6\:yy    \:divide\:dd\:x
        \else
           { \advance \:X  -0.5\:thickness
\ifdim \:X<\:LBorder \global\:LBorder\:X \fi
\advance \:X  \:thickness
\advance \:X  \:x
\ifdim \:X>\:RBorder \global\:RBorder\:X\fi }
           \let\:Yunitsy\:y  \let\:Xunitsx\:x  \let\:temp\:Vline
           \:dd0.6\:x
\ifdim \:x>\:mmp    \:divide\:dd\:yy   \fi
        \fi
        \advance\:dd  0.4\p@
\:ragged\:Cons\:dd\:ragged
        \:HVLn
      \fi }
   \advance\:X  \:x   \advance\:Y  \:y}\Define\:HVLn{
   \:xx\:AbsVal\:Xunitsx  \:divide\:xx\:ragged
   \advance\:xx  0.99\p@   \:K\:InCons\:xx \relax
   \ifnum \:K>\z@
      \divide\:Xunitsx \:K  \advance\:K \@ne
      \divide\:Yunitsy \:K
   \else \:K\@ne  \fi  \:NextLn}

\Define\:NextLn{
   \ifnum\:K=\z@  \let\:NextLn\relax
   \else  { \:temp }  \advance\:K  \m@ne
      \advance\:X  \:x   \advance\:Y  \:y
   \fi   \:NextLn}\newdimen\:ragged

\Define\Ragged(1){    \:ragged#1\p@  \:ragged0.1\:ragged }
\Ragged(7.5)
\Define\PaintUnderCurve(4){{
   \:Z\:Y   \def\:next{\Curve(#1,#2,#3,#4)}
   \MoveToLoc(#1)  \:d\:X   \MoveToLoc(#4)
   \advance\:d -\:X
   \ifdim \:AbsVal\:d<\:mmp  \def\:next{}
   \else
      \def\:CrvLnTo(##1,##2){
         \:x \:X   \:y  \:Y    \:X\:DJ  \:Y\:yyyy
         \:xx\:X   \:dddd \:Y  \:X\:x  \:Y\:y
         { \advance\:Y  \:dddd   \divide\:Y \tw@
           \advance\:Z -\:Y
           \advance\:Y   0.5\:Z
           \:dddd  \:AbsVal \:Z  \:d\z@
           \def\:CrvLnTo{\:LnTo}
           \:yy\:Y  \:dd\:dddd  \:ddd\:dddd
           \::paint  }}
    \fi \:next       }}
\Define\DoCurve(1){ \let\:StopCurve\:SlowCurve
  \def\:CMv(##1){  \:x\:X \:y\:Y   \:MvTo(##1)
      \advance\:x -\:X   \advance\:y -\:Y
      \:xxx \:x    \:yyy\:y}
   \:DoCurve{\Curve(#1)}
   \let\:StopCurve\:FastCurve}

\def\:DoCurve#1(#2)#3{{\XSaveUnits
   \def\:next{#1}    \:MvTo(#2,#2)
   \:x\:AbsVal\:X  \:y\:Y  \:ddd\z@  \:length
   \:Z\:d   \:divide\:Z{1.41421\p@}
   \edef\:tempa{\the\:DoDist}   \global\:DoDist\z@
   \def\:CrvLnTo(##1){ \MarkLoc(1^)    \:CMv(##1)
      { \MarkLoc(2^)   \:ddd\z@    \:length
\:dd  \:DoDist  \global\advance\:DoDist  \:d
\:ddd \:DoDist  \:divide\:ddd\:Z
\DoReg\:InCons\:ddd  \:Z\DoReg\:Z

      \ifdim \:Z>\:dd
         \advance\:Z -\:DoDist
\advance\:dd -\:DoDist
\:divide\:Z\:dd
\advance\:X \:Cons\:Z\:xxx
\advance\:Y \:Cons\:Z\:yyy   \:DoRot
         \def\:CrvLnTo{\:LnTo}
         \def\:OvalLn{\:Ln}  \XRecallUnits    #3 \fi}}
   \:next    \xdef\:DoDim{\:Cons\:DoDist}
   \global\:DoDist\:tempa     }
   \let\DoDim\:DoDim}

\newdimen\:DoDist\def\:DoRot{ \DSeg\RotateTo(1^,2^) }
\def\DoLine(#1,#2)(#3)#4{
   \MarkLoc($1)  \Move(#1,#2)
   \def\:next{   { \MarkLoc($2)
     \DSeg\RotateTo($1,$2)   \let\:DoRot\relax
     \edef\:RecallRagged{\the\:ragged} \MoveTo(#3,#3)
     \:x\:AbsVal\:X  \:y\:Y  \:ddd\z@  \:length
     \:ragged\:d   \divide\:ragged \tw@
     \DoCurve($1,$1,$2,$2)(#3)
        {\:ragged\:RecallRagged #4}  }
   \let\DoDim\:DoDim}  \:next }
\def\Table#1{\begingroup  \:DraCatCodes   \:multid#1
   \:DefineData#1}

\def\:DefineData#1#2{\endgroup
   \let\:temp~  \def~{\noexpand~}
   \edef#1{\noexpand\:DoPoly\expandafter\noexpand\csname :\string#1\endcsname}
   \expandafter\edef  \csname :\string#1\endcsname
       ##1{\noexpand\ifcase##1(#2)\noexpand\fi}
   \let~\:temp   \:DoNextPoly   \:DoNextPoly}

\def\:OR{\let\:or\or}  \:OR \catcode`\&13  \def&{)\noexpand\:or(}

\def\:TableData#1#2#3{\endgroup   \Table\:temp{#3}
   \:K\z@   \:J\z@    \def\:tempa(##1){\advance\:J \@ne }
   \:temp(0,999){\:tempa}    \let\:tempa&      \def\:temp{\def#1}
   \def&##1&{
      \ifnum  \:K<\:J
         \advance\:K \@ne
         \ifnum  \:K=\@ne   \def#1{#2(##1)}
         \else
            \:IIIexpandafter\:temp\expandafter{
               #1 & #2(##1) }
         \fi
      \else  \let\:next\relax \fi
      \:next}
   \let\:next&     &#3&&   \let&\:tempa }

\catcode`\&4  \def\:DoPoly#1(#2,#3)#4{
   \expandafter\let \csname :Back\the\:level\endcsname\:or
\expandafter\edef\csname :DoVars\the\:level\endcsname{
   \:DoB\the\:DoB}
\advance\:level  \@ne
   \:DoB#3  \advance\:DoB -#2
   \def\:PolyOr(##1){
      \ifnum  \:DoB=\z@  \:OR
      \else   #4(##1)   \advance\:DoB \m@ne    \fi}
   \:OR
   \def\:temp{\let\:or\:PolyOr #4}
   \:IIIexpandafter\:temp#1{#2}
   \advance\:level  \m@ne
\csname :DoVars\the\:level\endcsname
\def\:temp{\let\:or}
\expandafter\:temp\csname :Back\the\:level\endcsname }
\Define\PaintRect(2){\def\:next{{   \MarkLoc(^)
   \MoveToLoc(^) \Move( #1,0) \MarkLoc(^1) \Move(0,#2)
   \MarkLoc(^2)  \Move(-#1,0) \MarkLoc(^3)
   \PaintQuad(^,^1,^2,^3)        }}\:next}

\Define\PaintRectAt(4){\def\:next{{   \MoveTo(#1,#2)
   \:K#3   \advance\:K -#1
   \:DoB#4 \advance\:DoB -#2  \PaintRect(\:K,\:DoB)}}
   \:next}\def\::paint{
   \ifdim \:d<\:ragged      \advance\:xx -\:X
      \:yyy\:Y \:xxx\:dddd
      \advance\:Y \:yy  \divide\:Y \tw@
      \:average\:dddd\:dd\:ddd
      \def\:next{\:brush(\:xx,\z@)\:Y\:yyy\:dddd\:xxx}
   \else  \divide\:d \tw@
      \:average\:x\:X\:xx
      \:average\:y\:Y\:yy
      \:average\:dddd\:dd\:ddd
   \fi   \:next}\Define\:paint{{   \:AbsDif\:d\:xx\:x
   \ifdim \:d<\:mmp  \let\::paint\relax
   \else
      \ifdim  \:y >\:yyy   \:dd\:y   \:y \:yyy   \:yyy \:dd  \fi
      \ifdim  \:yy>\:yyyy  \:dd\:yy  \:yy\:yyyy  \:yyyy\:dd  \fi
      \:AbsDif\:dd\:yyyy\:yyy      \:AbsDif\:ddd\:yy\:y
      \ifdim \:dd<\:ddd  \:dd\:ddd  \fi
      \ifdim \:d >\:dd   \:d\:dd  \fi
      \advance\:d \:d      \:diff\:dd\:y\:yyy
      \:diff\:ddd\:yy\:yyyy
      \advance\:y    -0.5\:dd    \:abs\:dd
      \advance\:yy   -0.5\:ddd   \:abs\:ddd
      \:X\:x  \:Y\:y
      \:average\:dddd\:dd\:ddd
      \def\:next{ \:lpaint \:rpaint }
   \fi
   \::paint }}

\Define\:lpaint{ { \:xx\:x  \:yy \:y  \:ddd\:dddd  \::paint} }
\Define\:rpaint{ { \:X \:x  \:Y  \:y  \:dd \:dddd  \::paint} }
\Define\PaintQuad(4){\def\:next{{\Units(1pt,1pt)
   \MoveToLoc(#1)  \:x  \:X  \:y  \:Y
   \MoveToLoc(#2)  \:xx \:X  \:yy \:Y
   \MoveToLoc(#3)  \:xxx\:X  \:yyy\:Y
   \MoveToLoc(#4)  \:xxxx  \:X  \:yyyy \:Y

   \:paintQuad  }}\:next}

\def\:paintQuad{{
   \:SetVal\:a\:x\:y\:xx\:yy\:xxxx\:yyyy
\:SetVal\:b\:xx\:yy\:xxx\:yyy\:x\:y
\:SetVal\:c\:xxx\:yyy\:xx\:yy\:xxxx\:yyyy
\:SetVal\:cc\:xxxx\:yyyy\:xxx\:yyy\:x\:y
\def\:A{\:a} \def\:B{\:b} \def\:C{\:c} \def\:D{\:cc}
\:sort\:B\:A
\:sort\:C\:B  \:sort\:B\:A
\:sort\:D\:C  \:sort\:C\:B  \:sort\:B\:A
\let\:temp\relax
\:IsTriang\:A\:B
\:IsTriang\:B\:C
\:IsTriang\:C\:D
\:temp
   \:Quad\:A\:temp>   \:xxxx\:xx \:yyyy\:yy \:Z\:d
\:Quad\:D\:next<
   \:PrePaint(\:xx,\:yy,\:d,\:xxxx,\:yyyy,\:Z)
   \:temp  \:next }}
\Define\:PrePaint(6){
   \:x#1 \:y#2 \:yyy#3 \:xx#4 \:yy#5 \:yyyy#6
   \:paint }\def\:Quad#1#2#3{
         \:GetVal#1\:x\:y0
   \:GetVal#1\:xx\:yy1
   \:GetVal#1\:xxx\:yyy2
   \ifdim \:xx#3\:xxx
      \:ddd\:xx   \:xx\:xxx   \:xxx\:ddd
      \:ddd\:yy   \:yy\:yyy   \:yyy\:ddd
   \fi
          \def#2{}
   \:diff\:dd\:xxx\:xx
   \ifdim \ifdim\:AbsVal\:dd<\:mp  \:yy=\:yyy   \else  \z@>\z@  \fi

      \:d\:yy
   \else
      \ifdim \:AbsVal\:dd>\:mp
         \:diff\:dd\:xxx\:x   \:diff\:ddd\:yyy\:y
         \:divide\:ddd\:dd    \:diff\:dd\:xx\:xxx
         \:ddd\:Cons\:ddd\:dd   \advance\:ddd \:yyy
         \:d\:ddd
      \else \:d\:yyy \fi
      \edef#2{  \noexpand\:PrePaint
         (\the\:x ,\the\:y ,\the\:y,
         \the\:xx,\the\:yy,\the\:d)    }
   \fi}\def\:SetVal#1#2#3#4#5#6#7{
   \edef#1{(\the#2,\the#3,\the#4,\the#5,\the#6,\the#7)}}

\def\:sort#1#2{
   \ifdim \:IIIexpandafter\:field#1 <
          \:IIIexpandafter\:field#2
      \let\:temp#1  \let#1#2  \let#2\:temp
   \fi  }

\Define\:field(6){#1}

\def\:GetVal#1#2#3{
   \:IIIexpandafter\::GetVal #1#2#3}

\def\::GetVal(#1,#2,#3,#4,#5,#6)#7#8#9{
      \ifcase #9 #7#1   #8#2\or #7#3   #8#4\or #7#5   #8#6 \fi}
\def\:IsTriang#1#2{
   \ifdim \:IIIexpandafter\:field#1 =
          \:IIIexpandafter\:field#2
      \ifdim \:IIIexpandafter\:fieldB#1 =
             \:IIIexpandafter\:fieldB#2
         \def\:temp{ \:FixTria }
   \fi \fi  }

\def\:FixTria{
   \edef\:temp{\:IIIexpandafter\:FrsII\:B}
   \ifdim \:IIIexpandafter\:field\:A =
          \:IIIexpandafter\:field\:B
      \ifdim \:IIIexpandafter\:fieldB\:A =
             \:IIIexpandafter\:fieldB\:B
          \edef\:temp{\:IIIexpandafter\:FrsII\:C}
   \fi\fi
   \edef\:A{\:IIIexpandafter\:FrsII\:A}
   \edef\:D{\:IIIexpandafter\:FrsII\:D}
   \edef\:temp{
      \def\noexpand\:a{(\:A,\:temp,\:D)}
      \def\noexpand\:b{(\:temp,\:A,\:D)}
      \def\noexpand\:c{\noexpand\:b}
      \def\noexpand\:cc{(\:D,\:A,\:temp)}}
   \:temp
   \def\:A{\:a}  \def\:B{\:b}  \def\:C{\:c}  \def\:D{\:cc}  }

\Define\:fieldB(6){#2}
\Define\:FrsII(6){#1,#2}
\def\:IIIexpandafter{\expandafter\expandafter\expandafter}

\Define\DrawRect(2){ \Line( #1,0) \Line(0, #2)
                     \Line(-#1,0) \Line(0,-#2)}

\Define\DrawRectAt(4){{
   \MoveTo(#1,#2) \LineTo(#3,#2) \LineTo(#3,#4)
   \LineTo(#1,#4) \LineTo(#1,#2)}}
\def\du#1{ \ifx#1\:Xunits   \else\ifx#1\:Yunits
      \else\ifx#1\:Zunits   \else #1
      \fi \fi \fi}\Define\XSaveUnits{
   \expandafter\edef\csname XRecallUnits\the\:level\endcsname{
      \:StoreUnits}
    \advance\:level  \@ne}

\Define\XRecallUnits{
   \advance\:level \m@ne
   \csname XRecallUnits\the\:level \endcsname}

\Define\SaveUnits{
     }

\Define\:StoreUnits{       \:Xunits \the\:Xunits
   \:Yunits \the\:Yunits \:Zunits \the\:Zunits
   \:Xunitsx\the\:Xunitsx \:Xunitsy\the\:Xunitsy
\:Yunitsx\the\:Yunitsx \:Yunitsy\the\:Yunitsy }
\Define\:SearchDir{
   \ifdim \:x<\z@  \:x-\:x  \:y-\:y
   \edef\:tempA{\advance\:ddd  \ifdim \:y<\z@ - \fi\:CVXXX}
\else         \def\:tempA{} \fi
\ifdim \:y<\z@
   \edef\:tempA{\:ddd-\:ddd \advance\:ddd  \:CCCVX \:tempA}
   \:y-\:y
\fi
\ifdim \:y>\:x
   \:ddd\:y \:y\:x \:x\:ddd
   \edef\:tempA{\advance\:ddd  -\:XC \:ddd-\:ddd \:tempA}
\fi
   \:divide\:y\:x  \:d57.29578\:y
   \:ddd\:d       \:sqr\:y  \:K  \@ne
   \Do(1,30){\advance\:K \tw@
      \:d-\:Cons\:y\:d  \:dd\:d  \divide\:dd \:K
      \advance\:ddd \:dd  }
   \advance\:ddd -0.49182\:dd
   \:tempA }\Define\Curve(4){\def\:next{{  \XSaveUnits \Units(1pt,1pt)
   \MoveToLoc(#1) \:DI \:X  \:DK \:Y
   \MoveToLoc(#2) \:ddd \:X  \:Ez \:Y
   \MoveToLoc(#3) \:dd\:X  \:Ey\:Y
   \MoveToLoc(#4) \:DJ\:X  \:yyyy\:Y   \:Curve    }}\:next}
\Define\:FastCurve{   \:AbsDif\:d\:DI\:ddd  \:AbsDif\:dddd\:DK\:Ez
   \advance\:d \:dddd
   \ifnum \:d<\:ragged
      \:AbsDif\:d\:DJ\:dd \:AbsDif\:dddd\:yyyy\:Ey
      \advance\:d \:dddd                     \fi}
\let\:StopCurve\:FastCurve
\Define\:SlowCurve{
   \:AbsDif\:d\:DI\:DJ  \:AbsDif\:dddd\:DK\:yyyy
   \advance\:d \:dddd                     }
\Define\:Curve{   \:StopCurve
   \ifnum \:d<\:ragged   \:X\:DI \:Y\:DK
                    \:CrvLnTo(\:Cons\:DJ,\:Cons\:yyyy)
   \def\:SubCurves{}  \fi
   \:SubCurves}

\def\:CrvLnTo{\:LnTo}\Define\:SubCurves{
   \:average\:yy\:DI\:ddd     \:average\:Ex\:DK\:Ez
   \:average\:ddd\:ddd\:dd    \:average\:Ez\:Ez\:Ey
   \:average\:dd\:dd\:DJ  \:average\:Ey\:Ey\:yyyy
   \:average\:Zunits\:yy\:ddd    \:average\:Vdirection\:Ex\:Ez
   \:average\:ddd\:ddd\:dd    \:average\:Ez\:Ez\:Ey
   \:average\:DL\:Zunits\:ddd  \:average\:xxxx\:Vdirection\:Ez
   { \:ddd \:yy  \:Ez \:Ex   \:dd\:Zunits
     \:Ey\:Vdirection \:DJ\:DL \:yyyy\:xxxx \:Curve }
   \:DI \:DL \:DK \:xxxx               \:Curve   }
\def\MoveToCurve[#1]{
   \Define\:BiSect(3){\MoveToLoc(##1)
      \CSeg[#1]\Move(##1,##2)
   \MarkLoc(##3) }\:MvToCrv}

\Define\:MvToCrv(4){  \:BiSect(#1,#2,:a) \:BiSect(#2,#3,:b)
   \:BiSect(#3,#4,:c) \:BiSect(:a,:b,:A) \:BiSect(:b,:c,:B)
   \:BiSect(:A,:B,:Q)}
\Define\:OvalDir(3){   \:CosSin{#3\p@}
   \:Zunits\p@   \:d#2\:Zunits   \:x\:Cons\:d\:x
   \:Zunits\p@   \:d#1\:Zunits   \:y\:Cons\:d\:y
   \:SearchDir   }\Define\DrawOval   (2){ \DrawOvalArc(#1,#2)(0,\:cccvx) }
\Define\PaintOval  (2){ \PaintOvalArc(#1,#2)(0,\:cccvx) }
\Define\DrawCircle (1){ \DrawOval(#1,#1) }
\Define\PaintCircle(1){ \PaintOval(#1,#1) }
\def\DrawOvalArc(#1,#2)(#3,#4){{
       \:xxxx#4\p@  \advance\:xxxx -#3\p@
       \ifdim \:xxxx=\z@ \else
   \let\:SinOne\:SinB  \:OvalDir(#1,#2,#3)  \:DJ\:ddd
\:OvalDir(#1,#2,#4)  \:diff\:DI\:ddd\:DJ
\ifdim\:DI<\z@ \advance\:DI  \:CCCVX \fi
\ifdim \:xxxx<\:CCCVX \else \:DI\:CCCVX \fi
\:InitOval(#1,#2)  \:CosSin\:DJ
   \:xxxx\:x  \:yyyy\:y   \:xx\:X  \:yy\:Y
   \advance\:X \:Xx\:x  \advance\:X \:Yx\:y
   \advance\:Y \:Xy\:x  \advance\:Y \:Yy\:y
   \let\:Xunits\empty  \let\:Yunits\empty
   \Do(1,\:InCons\:DI){
      \:dd\:X  \:ddd\:Y  \:X\:xx  \:Y\:yy
      \:AdvOv\:xxx\:yyy\:xxxx\:yyyy
      \:X\:dd  \:Y\:ddd
      \advance\:xxx -\:X   \advance\:yyy -\:Y
      \:d\:AbsVal\:xxx
      \advance\:d  \:AbsVal\:yyy
      \ifdim \:d>\:ragged
         \:OvalLn(\:xxx,\:yyy)  \fi  }
   \:OvalDir(#1,#2,#4)    \:CosSin\:ddd
   \advance\:xx \:Xx\:x  \advance\:xx \:Yx\:y
   \advance\:yy \:Xy\:x  \advance\:yy \:Yy\:y
   \advance\:xx -\:X     \advance\:yy -\:Y
   \:OvalLn(\:xx,\:yy) \fi }}

\def\:OvalLn{\:Ln}\def\DoOvalArc(#1)(#2){   \:xx\:X  \:yy\:Y
   \def\:CMv(##1){  \:Mv(\:xxx,\:yyy)
      \:x\:xxx  \:y\:yyy}
   \:DoCurve{             \:X\:xx  \:Y\:yy
       \def\:DoRot{  \let\:Xunits\:XunitsReg
                     \let\:Yunits\:YunitsReg
                     \DSeg\RotateTo(1^,2^)     }
       \let\::OvalLn\:CrvLnTo
\Define\:OvalLn(2){ \:dd\:AbsVal####1
   \advance\:dd \:AbsVal####2 \:divide\:dd\:ragged
   \:J\:InCons\:dd  \advance\:J  \@ne
   \divide####1  \:J  \divide####2  \:J
   \Do(1,\:J){\::OvalLn(####1,####2)}}
       \DrawOvalArc(#1)(#2)}}
\def\NextTable{\begingroup  \:DraCatCodes \:NextTable}

\def\:NextTable#1{\endgroup
  \def\:DoNextPoly{#1\NextTable{}}}
\NextTable{}\def\:AdvOv#1#2#3#4{
   \:d\:CosOne#3 \advance\:d -\:SinOne#4
   #4\:CosOne#4    \advance#4     \:SinOne#3   #3\:d
   \divide#3 \:eight  \divide#4 \:eight
   #1\:X   #2\:Y
   \:d\:Xx#3  \advance\:d \:Yx#4  \advance#1  \:d
   \:d\:Xy#3  \advance\:d \:Yy#4  \advance#2  \:d  }

\def\:CosOne{7.99878}   \def\:SinB{0.13962}
\def\PaintOvalArc(#1,#2)(#3,#4){{ \ifdim #3\p@=#4\p@
                      \let\:next\relax  \else
   \:d\:AbsVal{#1\:Xunits} \advance\:d \:AbsVal{#2\:Yunits}
   \ifdim \:d<3\:ragged   \divide \:d \tw@   \PenSize(\:d)
\:Mv(-0.5\:d\du,0)  \:Ln(\:d\du,0)
   \else    \:InitOval(#1,#2)
      \MarkLoc(o$)   \RotateTo(#3) \MoveFToOval(#1,#2)
      \:Ex\:X  \:Ey\:Y   \edef\:FirstOvDir{\:Cons\:ddd\p@}
      \MoveToLoc(o$) \RotateTo(#4) \MoveFToOval(#1,#2)
      \:Ez \:X  \:Vdirection \:Y   \edef\:LastOvDir{\:Cons\:ddd\p@}
      \MoveToLoc(o$)
      \if:rotated
      \:Zunits\p@     \:Zunits#1\:Zunits
      \:xx\:Cons\:Zunits\:Xunitsx
      \:Zunits\p@     \:Zunits#2\:Zunits
      \:yy\:Cons\:Zunits\:Yunitsx
      \:x\:xx  \:y\:yy  \:ddd\z@  \:length
      \:ddd\:d  \:divide\:xx\:ddd   \:divide\:yy\:ddd
\else \:xx\p@  \:yy\z@   \fi
      \:AbsDif\:d{#3\p@}{#4\p@}
      \ifdim \:d>359\p@ \:Ez-\:Xx\:xx  \advance\:Ez -\:Yx\:yy
\:Vdirection-\:Xy\:xx  \advance\:Vdirection -\:Yy\:yy
\advance\:Ez \:X  \advance\:Vdirection \:Y
\:setpaint\:PaintOvOv<>
      \else
         \:x\:xx \:y\:yy \:SearchDir
\:xxx\:FirstOvDir  \advance\:xxx -\:ddd
\ifdim \:xxx<\z@    \advance\:xxx \:CCCVX  \fi
\:yyy\:LastOvDir   \advance\:yyy -\:ddd
\ifdim \:yyy<\z@    \advance\:yyy \:CCCVX  \fi
\:J\z@
\ifdim \:xxx<\:yyy  \ifdim        \:yyy<\:CVXXX
            \:Pntovln\:FirstOvDir\:FirstOvDir\:LastOvDir
\else \ifdim  \:xxx>\:CVXXX
            \:Pntovln\:FirstOvDir\:FirstOvDir\:LastOvDir
\else \:yyy-\:yyy  \advance\:yyy \:CCCVX
      \ifdim  \:xxx<\:yyy
            \:setpaint\:PntLeftOvOv><  \:PntMovln\:FirstOvDir\:xx
      \else \:FxLx  \:setpaint\:PntLeftOvOv><  \:Pntmovln\:LastOvDir
\fi  \fi  \fi
\else               \ifdim        \:xxx<\:CVXXX
            { \:setpaint\:PaintOvOv<> }  \:FxLx  \:setpaint\:PntLeftOvOv><
\:Usrch  \:PaintMidOvLn\:LastOvDir\:FirstOvDir
\else \ifdim  \:yyy>\:CVXXX
            {  \:FxLx \:setpaint\:PaintOvOv<>  }  \:setpaint\:PntLeftOvOv><
 \:Dsrch  \:PaintMidOvLn\:LastOvDir\:FirstOvDir
\else \:xxx-\:xxx  \advance\:xxx \:CCCVX
      \ifdim  \:yyy<\:xxx
            \:setpaint\:PaintOvOv<>    \:PntMovln\:FirstOvDir\:xx
      \else
            \:FxLx  \:setpaint\:PaintOvOv<>   \:Pntmovln\:LastOvDir
\fi  \fi  \fi  \fi
   \fi \fi \fi}}\def\:setpaint#1#2#3{{\aftergroup#1
\:d\:Xx\:xx   \advance\:d \:Yx\:yy
\ifdim \:d<\z@  \aftergroup#3
\else           \aftergroup#2  \fi}}
\def\:FxLx{\:d\:Ex  \:Ex\:Ez  \:Ez\:d}

\def\:Pntovln#1{
   \let\:SinOne\:SinB     \:CosSin#1
   \:xxxx\:Ex  \:yyyy\:Ey  \:Z\z@
   \:PaintOvLn}

\def\:PntMovln{
   \:Dsrch  \:FxLx   \:xx\:ddd  \:PaintOvLn}

\def\:Pntmovln{
   \:Usrch  \:FxLx  \:xx\:ddd   \:PaintOvLn\:xx}
\def\:PaintMidOvLn#1#2{
   \:FxLx   \:xx\:ddd  \:PaintOvLn#1#2
   \:xxx\:Ez  \:yyy\:Vdirection
   \:ddd\:Xy\:x  \advance\:ddd \:Yy\:y  \advance\:ddd \:Y
   \:PaintSlice}

\def\:PntLeftOvOv{
   \:xx-\:xx  \:yy-\:yy  \:PaintOvOv}

\def\:Dsrch{      \def\:SinOne{-\:SinB}
   \:xx\:x  \:yy\:y   \:SearchDir
   \:x\:xx  \:y\:yy
   \:d\:Ey  \:Ey\:Vdirection  \:Vdirection\:d }

\def\:Usrch{
   \let\:SinOne\:SinB
   \:x\:xx  \:y\:yy   \:SearchDir
   \:x\:xx  \:y\:yy}\def\:PaintOvLn#1#2{
   \:diff\:dd\:Ey\:Vdirection  \:diff\:ddd\:Ex\:Ez
   \ifdim \:AbsVal\:ddd>\:mp
      \:divide\:dd\:ddd      \:d#2
      \advance\:d -#1
      \ifnum \:d<\z@  \advance\:d \:CCCVX  \fi
      \:DoB\:InCons\:d  \let\:next\:PntDo  \:next
   \fi}

\def\:PntDo{
   \ifnum\:DoB=\z@ \let\:next\relax
   \else
      \::AdvOv\:x\:y
      \ifdim \:d>\:ragged
         \:ddd\:xxx  \advance\:ddd -\:Ez
         \:ddd\:Cons\:dd\:ddd
         \advance\:ddd \:Vdirection  \:PaintSlice
      \fi
      \advance\:DoB \m@ne
   \fi  \:next}\Define\:PaintSlice{   \:AbsDif\:dddd\:yyy\:ddd
   \advance\:yyy \:ddd  \divide\:yyy \tw@
   { \advance\:dddd \:Z  \divide\:dddd \tw@
      \:X\:xxxx  \:Y\:yyy
      \:xx\:xxx  \advance\:xxx -\:xxxx
      \:d\z@
      \:yy\:Y  \:dd\:dddd  \:ddd\:dddd
      \::paint     }
   \:xxxx\:xxx  \:yyyy\:yyy  \:Z\:dddd    \:J\z@}
\def\::AdvOv#1#2{  \:AdvOv\:xxx\:yyy#1#2
   \:AbsDif\:d\:xxxx\:xxx       \advance\:J \@ne
   \ifnum \:J=\sixt@@n   \multiply\:d \@cclvi
   \fi }\def\:PaintOvOv#1{   \def\:hdir{#1}
   \:xxx\:Xx\:xx  \advance\:xxx \:Yx\:yy
   \:yyy\:Xy\:xx  \advance\:yyy \:Yy\:yy
   \advance\:xxx \:X   \advance\:yyy \:Y
   \:Z\z@  \:xxxx\:xxx  \:yyyy\:yyy
   \:x\:xx   \:y\:yy   \:DoB\z@   \:J\z@
   \let\:next\:scanOvOv  \:next }
\Define\:scanOvOv{  \advance\:DoB \@ne
   \advance\:J \@ne
   \:AbsDif\:d\:xxx\:Ez
   \ifdim \ifdim      \:xxx\:hdir\:Ez    -
          \else\ifdim \thr@@\:d<\:ragged -
          \else\ifnum \:DoB>358            -
          \fi \fi \fi                     \p@<\z@
      \:ddd\:yyyy  \advance\:ddd -0.5\:dddd
      \:yyy\:yyyy   \advance\:yyy   0.5\:dddd
      \:xxx\:Ez \:PaintSlice
      \let\:next\relax
   \else
     \:d\:xx \advance\:d -\:x
     \:xxx\:yy \advance\:xxx -\:y
     \ifdim -\:Xx\:d\:hdir\:Yx\:xxx
        \let\:SinOne\:SinB
        \::AdvOv\:xx\:yy
        \ifdim  \ifdim       \:d >\:ragged -
                \else \ifnum \:J>\sixt@@n -
                \fi\fi                      \p@<\z@
           \:J\z@  \def\:SinOne{-\:SinB}
           \:AdvOv\:dd\:ddd\:x\:y  \:PaintSlice  \fi
     \else
        \def\:SinOne{-\:SinB}
        \::AdvOv\:x\:y
        \ifdim  \ifdim       \:d >\:ragged -
                \else \ifnum \:J>\sixt@@n -
                \fi\fi                      \p@<\z@
           \:J\z@   \let\:SinOne\:SinB
           \:AdvOv\:dd\:ddd\:xx\:yy  \:PaintSlice  \fi
   \fi \fi
   \:next}\Define\SetBrush{\:Opt[]\:SetBrush{}}

\def\:SetBrush[#1](#2,#3)#4{    \def\:temp{#4}
   \ifx \:temp\empty
      \def\:brush{   \let\:Xunits\empty \let\:Yunits\empty
                     \:thickness\:dddd  \:Ln  }
   \else       \def\:BruShape{#4}
      \:dd#2\:Xunits      \:ddd#3\:Yunits
      \edef\:Grd{ \:dd\the\:dd  \:ddd\the\:ddd }
      \MarkLoc($$)  \def\:temp{#1}
\ifx \:temp\empty  \:X\z@  \:Y\z@
\else  \MoveTo(#1) \fi
\edef\:BrOrg{ \:x\the\:X  \:y\the\:Y }
\MoveToLoc($$)
      \def\:brush(##1,##2){ \::brush }  \fi  }

\SetBrush(,){}\def\::brush{{  \SetBrush(,){}
\let\:Xunits\:XunitsReg   \let\:Yunits\:YunitsReg
\advance\:Y -0.5\:dddd   \:yy\:Y
\advance\:yy \:dddd
\:BrOrg  \:Grd  \advance\:xx  \:X
\ifdim \:xx<\:X  \:d\:X \:X\:xx \:xx\:d  \fi
   \:GridPt\:X\:x\:dd
   \:GridPt\:Y\:y\:ddd            \:x\:X
   \:DoBrush       }}\Define\:DoBrush{
   \ifdim       \:Y>\:yy  \let\:DoBrush\relax
   \else \ifdim \:X>\:xx \advance\:Y \:ddd     \:X\:x
   \else { \:BruShape }  \advance\:X \:dd
   \fi  \fi  \:DoBrush  }\def\:GridPt#1#2#3{   \:xxxx#1
   \advance#1 -#2  \:divide#1#3
   #1\:InCons#1#3  \advance#1 #2
   \ifdim #1=\:xxxx
   \else  \ifdim \:xxxx>#2 \advance#1 #3 \fi \fi  }
\newcount\:IntId    \newcount\:IntCount
\newcount\:DecId    \newcount\:DecCount

\Define\:NewCount{\alloc@ 0\count \countdef \insc@unt }
\Define\:NewDimen{\alloc@ 1\dimen \dimendef \insc@unt }

\def\:NewVar#1#2#3#4#5{ \:multid#1
   \def\:temp{    \csname \string#4\the#4\endcsname\z@
      \edef#1{\noexpand#3  \csname \string#4\the#4\endcsname}}
   \def\:next{  \global#5#4   \expandafter
       #2   \csname \string#4\the#4\endcsname  \:temp }
   \advance#4  \@ne
   \ifnum #4 > #5 \else \def\:next{\:temp}  \fi   \:next}

\def\IntVar#1{\:NewVar#1\:NewCount\:IntOp\:IntId\:IntCount}
\def\DecVar#1{\:NewVar#1\:NewDimen\:DecOp\:DecId\:DecCount}

\DecVar\Q  \DecVar\R   \DecVar\T
\IntVar\I  \IntVar\J   \IntVar\K
\def\WriteVal#1{\immediate\write\sixt@@n{...\string#1=#1;}}

\newdimen\:X   \newdimen\:Y
\newdimen\:x   \newdimen\:y   \newdimen\:d
\newdimen\:xx  \newdimen\:yy  \newdimen\:dd
\newdimen\:xxx \newdimen\:yyy \newdimen\:ddd
\newdimen\:xxxx\newdimen\:yyyy\newdimen\:dddd
\newcount\:J  \newcount\:K
\newdimen\:DI   \newdimen\:DJ
\newdimen\:DK   \newdimen\:DL   \newtoks\:t
\def\:IntFromPt#1#2{
   \:d#2\relax
   \advance\:d  \ifdim\:d<-0.5\p@-\fi  0.5\p@
   #1\:d    \divide#1  65536\relax}
\def\:temp{\catcode`\p12  \catcode`\t12}
\def\:Cons{\catcode`\p11  \catcode`\t11}
\:temp  \def\:Frac#1pt{#1}
        \def\:rnd#1.#2pt{#1}  \:Cons

\def\:Cons#1{\expandafter\:Frac\the#1}
\def\:sqr#1{#1\expandafter\:Frac\the#1#1}
\def\:InCons#1{\expandafter\:rnd\the#1}\def\:Val#1{#1;}
\let\Val\:Val\def\:IntOp#1#2{\csname :Op#2\endcsname#1}
\let\:SvIntOp\:IntOp
\def\:PreIntOp{\let\:IntOp\empty
   \let\Val\empty}
\def\:PostIntOp{\let\:IntOp\:SvIntOp
   \let\Val\:Val}

\expandafter\def\csname :Op;\endcsname#1{ \the#1}
\expandafter\def\csname :Op=\endcsname#1#2;{
   \:PreIntOp#1#2\:PostIntOp}
\expandafter\def\csname :Op+\endcsname#1#2;{
   \:PreIntOp\advance #1  #2\:PostIntOp}
\expandafter\def\csname :Op-\endcsname#1#2;{
   \:PreIntOp\advance #1  -#2\:PostIntOp}
\expandafter\def\csname :Op/\endcsname#1#2;{
   \:PreIntOp\divide#1   #2\:PostIntOp}
\expandafter\def\csname :Op*\endcsname#1#2;{
   \:PreIntOp\multiply#1   #2\:PostIntOp}
\def\:DecOp#1#2{ \csname :xOp#2\endcsname#1}
\let\:SvDecOp\:DecOp
\def\:PreDecOp{\let\:IntOp\the \def\:DecOp{\:Cons}
   \let\Val\empty   \let\:du\empty}
\def\:PostDecOp{\let\:IntOp\:SvIntOp \let\Val\:Val
   \let\:DecOp\:SvDecOp  \let\:du\::du  }

\def\::du#1{\p@
   \ifx#1\p@ \let\:temp\relax
   \else     \def\:temp{\du{#1}}
   \fi\:temp}                    \:PostDecOp

\expandafter\def\csname :xOp;\endcsname#1{ \:Cons#1}
\expandafter\def\csname :Op[\endcsname#1#2];{
   \:PreDecOp \:dd#2\p@  \:IntFromPt#1\:dd
                                 \:PostDecOp  }
\expandafter\def\csname :xOp=\endcsname#1#2;{
   \:PreDecOp#1#2\p@\:PostDecOp               }
\expandafter\def\csname :xOp(\endcsname#1#2){
   \:PreDecOp#1#2\p@\:PostDecOp               }
\expandafter\def\csname :xOp+\endcsname#1#2;{
   \:PreDecOp\advance #1  #2\p@\:PostDecOp  }
\expandafter\def\csname :xOp-\endcsname#1#2;{
   \:PreDecOp\advance #1  -#2\p@\:PostDecOp }
\expandafter\def\csname :xOp*\endcsname#1#2;{
   \:PreDecOp#1 #2#1\:PostDecOp              }
\expandafter\def\csname :xOp/\endcsname#1#2;{
   \:PreDecOp  \:divide#1{#2\p@}  \:PostDecOp }
\let\IF\ifnum

\def\EqText(#1,#2){
   \z@=\z@ \fi  \def\:temp{#1}
                \def\:next{#2}    \ifx \:temp\:next }

\def\:IfInt#1(#2,#3){ \z@=\z@ \fi
   \:IntOp\:K=#2;  \:IntOp\:J=#3; \ifnum  \:K#1\:J }

\def\:IfDim#1(#2,#3){ \z@=\z@ \fi
   \:DecOp\:d=#2;  \:DecOp\:dd=#3; \ifdim  \:d#1\:dd }

 \def\Do(#1,#2)#3{
   \expandafter\let
   \csname :Back\the\:level\endcsname\:Do
\expandafter\edef\csname :DoVars\the\:level\endcsname{
   \DoReg\the\DoReg \:DoB\the\:DoB}
\advance\:level  \@ne
   \DoReg#1  \:DoB#2  \relax
   \ifnum \DoReg<\:DoB
      \def\:Do{\ifnum \DoReg>\:DoB
                  \let\:Do\relax
               \else  #3\advance\DoReg  \@ne \fi
               \:Do}
   \else
      \def\:Do{\ifnum \DoReg<\:DoB
                  \let\:Do\relax
               \else  #3\advance\DoReg \m@ne  \fi
               \:Do}
   \fi  \def\:nextdo{ \:Do \advance\:level  \m@ne
\csname :DoVars\the\:level\endcsname
\def\:temp{\let\:Do}
\expandafter\:temp\csname
   :Back\the\:level\endcsname  } \:nextdo}

\let\:Do\relax

\newcount\DoReg   \let\:DoReg\DoReg       \newcount\:DoB
 \newcount\:level\def\::divide#1{   \:DI\:DK   \:dddd\:DL
   \advance\:DI -\:Cons\:dddd#1
   \:IntFromPt\:J\:dddd  \advance\:dddd -\:J\p@
   \multiply\:J  \@M   \:IntFromPt\:K{\@M\:dddd}
   \advance\:J \:K     \:dddd\@M\p@
   \divide\:dddd \:J   \advance#1 \:Cons\:DI\:dddd  }

\def\:divide#1#2{   \:DK#1   \:DL#2   #1\z@
   \::divide#1  \::divide#1  \::divide#1
   \::divide#1  \::divide#1  }
\def\:Sqrt#1{ \ifdim #1<\:mmp   #1\z@  \else
   \:dd#1   \divide\:dd \tw@
   \def\::Sqrt{  \:ddd#1
      \:divide\:ddd\:dd      \:AbsDif\:d\:dd\:ddd
      \advance\:dd \:ddd   \divide\:dd \tw@
      \ifdim  \:d < \:mmmp
         \let\::Sqrt\relax  \fi
      \::Sqrt}
   \::Sqrt   #1\:dd   \fi }\Define\:length{
   \:dd \:AbsVal \:ddd
   \:abs\:x  \:abs\:y
   \ifdim \:dd<\:x \:dd\:x \fi
   \ifdim \:dd<\:y \:dd\:y \fi
   \ifdim \:dd>\:mp
      \:divide\:x\:dd     \:sqr\:x
      \:divide\:y\:dd     \:sqr\:y
      \:divide\:ddd\:dd   \:sqr\:ddd
      \advance\:x  \:y  \advance\:x \:ddd
      \:y\:dd  \:Sqrt\:x  \:d\:Cons\:x\:y
   \else \:d\:dd \fi}\Define\:distance(2){\MarkLoc(@^)
   \MoveToLoc(#1)  \:x \:X  \:y \:Y  \:ddd\:Z
   \MoveToLoc(#2)
   \advance\:x  -\:X   \advance\:y  -\:Y
   \advance\:ddd  -\if:IIID \:Z \else \:ddd \fi
   \:length  \MoveToLoc(@^)}\def\:NormalizeDeg#1{
   \:DL#1   \:K\:InCons\:DL
   \divide\:K  \:cccvx   \multiply\:K  \:cccvx
   \advance #1 -\:K\p@
   \ifdim #1<\z@ \advance #1  \:CCCVX \fi
   \ifdim #1=\z@
      \ifdim\:DL=\z@ \else
         \advance #1  \:CCCVX \fi \fi }\def\:CosSin#1{ \:DK#1
   \:NormalizeDeg\:DK \def\:tempA{}
\ifdim \:CVXXX<\:DK
   \def\:tempA{\:y-\:y}
   \advance\:DK -\:CCCVX  \:DK-\:DK   \fi
\ifdim \:XC<\:DK
   \edef\:tempA{\:x-\:x \:tempA}
   \advance\:DK -\:CVXXX  \:DK-\:DK   \fi
\ifdim 45\p@<\:DK
   \edef\:tempA{\:d\:x \:x\:y \:y\:d \:tempA}
   \advance\:DK -\:XC   \:DK-\:DK   \fi
   \:x\p@   \:y0.01745\:DK   \:d\:y   \:K\@ne
   \edef\:next{\advance\:K \@ne
      \:sqr\:d  \divide\:d \:K  \advance}
   \:next \:x -\:d   \:next \:y -\:d
   \:next \:x  \:d   \:next \:y  \:d
   \:next \:x -\:d   \:next \:y -\:d
   \:next \:x  \:d   \:next \:y  \:d
   \:tempA   }   \Define\:rInitOval(2){
   \:Zunits\p@   \:dd#1\:Zunits
   \:d\:Cons\:dd\:Xunitsx   \edef\:Xx{\:Cons\:d}
   \:d\:Cons\:dd\:Xunitsy   \edef\:Xy{\:Cons\:d}
   \:dd#2\:Zunits
   \:d\:Cons\:dd\:Yunitsx   \edef\:Yx{\:Cons\:d}
   \:d\:Cons\:dd\:Yunitsy   \edef\:Yy{\:Cons\:d}}
\Define\:xyInitOval(2){
   \:d#1\:Xunits   \edef\:Xx{\:Cons\:d} \def\:Xy{0}
   \:d#2\:Yunits   \edef\:Yy{\:Cons\:d} \def\:Yx{0} }
 \def\:FigSize#1#2#3{
   \:x\:LBorder   \:y\:RBorder   \:d\:TeXLoc
   {\Object\:temp{#3}
    \setbox\:box\hbox{ \:temp
       \multiply\:x by \tw@  \multiply\:y by \tw@
       \xdef\:FSize{ \noexpand#1=\:Cons\:x;
                       \noexpand#2=\:Cons\:y;}}}
   \global\:LBorder\:x   \global\:RBorder\:y
   \global\:TeXLoc \:d
   \:FSize}
\expandafter\let \csname 0:Ln \endcsname\:Ln
\expandafter\def\csname 1:Ln \endcsname{
      \advance\:x -\:X  \advance\:y -\:Y
      \csname 0:Ln \endcsname(\:x,\:y)  }

\newcount\:ClipLevel   \:ClipLevel\@ne

\Define\Clip{\futurelet\:next\:Clip}

\Define\:Clip{
   \ifx \:next[  \expandafter\:DefClipOut
        \else    \expandafter\:DefClip     \fi }

\def\:DefClipOut[#1]{ \:DefClip(#1) }
\Define\:DefClip(1){  \def\:temp{#1}
   \ifx \:temp\empty
      \:ClipLevel\@ne     \def\:next{\let\:Ln}
      \expandafter\:next\csname 0:Ln \endcsname
   \else  \def\:temp{\::DefClip(#1)}  \fi  \:temp }
\Define\::DefClip(1){  \MarkLoc(^)
   \:x\:X  \:y \:Y  \Move(#1)
   \ifdim\:x>\:X  \:dd\:X  \:X\:x  \:x\:dd  \fi
   \ifdim\:y>\:Y  \:dd\:Y  \:Y\:y  \:y\:dd  \fi
   \advance\:ClipLevel  \@ne
   \expandafter\edef\csname \the
      \:ClipLevel :Ln \endcsname{
          \:xxx \the\:x   \:yyy \the\:y
          \:xxxx\the\:X  \:yyyy\the\:Y
          \ifx \:next[   \noexpand\:ClipOut
          \else          \noexpand\:ClipIn     \fi }
   \let\:Ln\:ClipLn   \MoveToLoc(^)  }\def\:ClipLn(#1,#2){
   \:x#1\:Xunits \:y#2\:Yunits
   {  \let\:Xunits\empty  \let\:Yunits\empty
   \advance\:x \:X   \advance\:y \:Y
   \ifdim \:x<\:X  \:dd\:X \:X\:x \:x\:dd
                   \:dd\:Y \:Y\:y \:y\:dd  \fi
   \:diff\:dd\:X\:x   \:diff\:ddd\:Y\:y
\:Z \:AbsVal \:dd
\advance\:Z  \:AbsVal\:ddd
\ifdim \:Z>\sixt@@n\p@
   \divide\:dd   128
   \divide\:ddd  128  \fi
\:Z\:Cons\:y\:dd
\advance\:Z -\:Cons\:x\:ddd
\ifdim \:dd<\z@
  \:dd-\:dd  \:ddd-\:ddd  \:Z-\:Z
\fi
   \csname \the \:ClipLevel :Ln \endcsname      }
   \advance\:X \:x   \advance\:Y \:y  }\Define\:ClipIn{
   \def\:next{\let\:next}
   \expandafter\:next\csname \the
      \:ClipLevel :Ln \endcsname
   \advance\:ClipLevel \m@ne
   { \:ClipLeft\:xxxx \:ClipDown\:yyy \:ClipUp\:yyyy \:next }
   { \:ClipRight\:xxx \:ClipDown\:yyy \:ClipUp\:yyyy \:next  }
   { \:ClipDown\:yyyy \:next }
     \:ClipUp\:yyy    \:next }

\Define\:ClipOut{
   \:ClipLeft\:xxx      \:ClipRight\:xxxx
   \:ClipUp  \:yyyy     \:ClipDown \:yyy
   \def\:next{\let\:next}
   \expandafter\:next
      \csname \the\:ClipLevel :Ln \endcsname
   \advance\:ClipLevel \m@ne      \:next  }\def\:ClipLeft#1{
   \ifdim       \:x<#1  \:KilledLine
   \else \ifdim \:X<#1  \:X#1
      \ifdim \:dd>\:mmmp
         \:Y\:Cons\:ddd\:X  \advance\:Y \:Z
         \:divide\:Y\:dd
   \fi \fi \fi     \:CondKilLn  }

\def\:ClipRight#1{
   \ifdim       \:X>#1  \:KilledLine
   \else \ifdim \:x>#1  \:x#1
      \ifdim \:dd>\:mmmp
         \:y\:Cons\:ddd\:x  \advance\:y \:Z
         \:divide\:y\:dd
   \fi \fi \fi    \:CondKilLn }

\Define\:CondKilLn{
   \:d\:x  \advance\:d -\:X
   \ifdim \:d<\z@  \:d-\:d \fi
   \ifdim \:y<\:Y  \advance\:d \:Y \advance\:d -\:y
   \else           \advance\:d \:y \advance\:d -\:Y  \fi
   \ifdim  \:d<\:mmp   \:KilledLine \fi
   \ifdim  \:thickness=\z@ \:KilledLine \fi}

\Define\:KilledLine{
   \let\:ClipLeft\:gobble  \let\:ClipRight\:gobble
   \let\:ClipUp  \:gobble  \let\:ClipDown \:gobble
   \let\:next\relax
   \expandafter\def\csname 1:Ln \endcsname{}}\def\:ClipUp#1{
   \:AbsDif\:d\:y\:Y
   \ifdim  \:d<\:ragged
      \advance\:y  0.5\:thickness
\advance\:Y -0.5\:thickness
\ifdim       \:Y>#1  \:KilledLine
\else \ifdim \:y>#1
   \:thickness#1 \advance\:thickness -\:Y
   \advance\:Y  0.5\:thickness  \:y\:Y
\else
   \advance\:Y  0.5\:thickness
   \advance\:y -0.5\:thickness
\fi  \fi
\:dd\p@   \:ddd\z@
\def\:temp{  \:Z\:Y }  \:temp
   \else \let\:temp\relax
          \ifdim \ifdim\:Y<\:y\:Y\else\:y\fi >#1  \:KilledLine
   \else  \ifdim  \::ClipUp#1\:X\:Y
   \else  \ifdim  \::ClipUp#1\:x\:y
   \fi \fi \fi \fi   \:CondKilLn  }

\def\::ClipUp#1#2#3{
#3>#1   #3#1
\ifdim \:AbsVal\:ddd>\:mmmp
   #2\:Cons\:dd#3  \advance#2 -\:Z
   \:divide#2\:ddd
\fi  }\def\:ClipDown#1{   \:Ex2#1
   \:Flip\:y  \:Flip\:Y  \:ClipUp#1
   \:Flip\:y  \:Flip\:Y  \:temp }

\def\:Flip#1{   #1-#1  \advance#1 \:Ex  }
\Define\:SetDrawWidth{
   \hsize\:RBorder      \advance\hsize -\:LBorder
   \leftskip -\:LBorder \rightskip\z@}
\newdimen\:thickness   \:thickness0.75\p@

\Define\PenSize(1){\:thickness#1\relax}
\Define\:Draw{   \ifvmode \noindent\hfil\fi
   \global\:TeXLoc\z@
   \vbox\bgroup
           \begingroup
   \def\EndDraw{
           \endgroup   \:SetDrawWidth
         \egroup}
   \:DraCatCodes      \parindent\z@    \everypar{}
\leftskip\z@     \rightskip\z@    \boxmaxdepth\maxdimen
\let\FigSize\:FigSize
\def\Draw{\:wrn1{}} \:CommonIID   \:InDraw }
\Define\:CommonIID{\def\RotateTo{\:RotateTo}

\def\MoveF{\:MvF}\def\LineToLoc{\:LnToLoc}}\newdimen\:Xunits   \:Xunits\p@
\newdimen\:Yunits   \:Yunits\p@
\Define\:InDraw{\:Opt()\::InDraw{\:Xunits,\:Yunits}}
\Define\::InDraw(2){   \:X\z@   \:Y\z@
   \:Xunits#1 \relax  \:Yunits#2  \:AdjRunits
   \global\:LBorder\z@    \global\:RBorder\z@
   \:loadIID   \leavevmode }

\Define\:Resize(2){
   \:Xunits#1\:Xunits  \relax
   \:Yunits#2\:Yunits  \:AdjRunits  }
\Define\:Units(2){
   \:Xunits#1 \relax   \:Yunits#2    \:AdjRunits   }
\Define\:loadIID{
   \def\LineAt{\:DLn}
   \def\LineTo{\:LnTo}
   \def\MoveTo{\:MvTo}
   \def\Line{\:Ln}
   \def\Move{\:Mv}
   \def\MoveF{\:MvF}
   \:rotatedfalse\def\:InitOval{\:xyInitOval}
   
   \def\Units{\:Units}}

\def\DrawOn{\def\Draw{\:Draw}}                      \DrawOn
\def\DrawOff{\def\Draw{\begingroup \:J\@cclv
                       \:NoDrawSpecials \:NoDraw}}
\catcode`\/0 \catcode`\\11
/def/:NoDraw#1\EndDraw{/endgroup}
/catcode`/\0 /catcode`//12

\def\:NoDrawSpecials{\catcode\:J11
  \ifnum \:J=\z@
     \let \:NoDrawSpecials\relax \fi
  \advance\:J  \m@ne \:NoDrawSpecials}
\let\:XunitsReg\:Xunits   \let\:YunitsReg\:Yunits  \Define\EntryExit(4){
   \edef\:InOut##1{
      \noexpand\ifcase ##1\space
         #1\noexpand\or #2\noexpand\or
         #3\noexpand\or #4\noexpand\fi}}

\EntryExit(0,0,0,0)\Define\:DrawBox{   \:x0.5\wd\:box
   \:y\ht\:box \advance\:y  \dp\:box
   \divide\:y \tw@
   \advance \:X -\:InOut0\:x
   \advance \:Y -\:InOut1\:y
   { \advance \:X -\:x   \advance \:Y -\:y
     { \ifdim\:X<\:LBorder \global\:LBorder\:X\fi
\advance\:X \wd\:box
\ifdim\:X>\:RBorder \global\:RBorder\:X\fi }
     \advance \:Y  \dp\:box
     {\:d\:X \advance\:d \wd\:box
 \advance\:X -\:TeXLoc   \global\:TeXLoc\:d
 \vrule width\:X depth\z@ height\z@
 \raise \:Y \box\:box} }
   \edef\MoveToExit(##1,##2){
   \:X\the\:X   \:Y\the\:Y
   \:x\the\:x   \:y\the\:y
   \advance\:X  ##1\:x
   \advance\:Y  ##2\:y}
   \advance \:X  \:InOut2\:x
   \advance \:Y  \:InOut3\:y }
 \Define\ThreeDim{\:Opt[]\:ThreeDim{\p@}}

\def\:ThreeDim[#1](#2){\::ThreeDim[#1](#2,,)}

\def\::ThreeDim[#1](#2,#3,#4,#5){ \bgroup\begingroup
   \def\EndThreeDim{          \endgroup\egroup}
   \:IIIDtrue   \:Zunits#1  \:Z\z@
   \def\:temp{#4}
   \ifx \:temp\empty   \:CosSin{#3\p@}  \:divide\:x\:y       \:Ey\:x
\:CosSin{#2\p@}  \:Ex\:Cons\:Ey\:x  \:Ey\:Cons\:Ey\:y
\let\:project\:projectPar

   \else               \:Ex#2\:Xunits \:Ey#3\:Yunits \:Ez#4\:Zunits
\let\:project\:projectPer
 \fi
   \def\LineAt{\:tDLn}
\def\LineTo{\:tLnTo}
\def\MoveTo{\:tMvTo}
\def\Line{\:tLn}
\def\Move{\:tMv}

\def\Units{\:tUnits}\def\RotateTo{\:tRotateTo}

\def\MoveF{\:tMvF}
\:Vdirection\z@ \def\LineToLoc{\:tLnToLoc} }  \Define\:projectPer{
   \:diff\:x\:X\:Ex   \:diff\:y\:Y\:Ey
   \:diff\:xxxx\:Z\:Ez
   \:divide\:x\:xxxx      \:divide\:y\:xxxx
   \:x-\:Cons\:Ez\:x     \:y-\:Cons\:Ez\:y
   \advance\:x  \:Ex    \advance\:y \:Ey}
\Define\:projectPar{
   \:x\:Cons\:Ex\:Z  \advance\:x \:X
   \:y\:Cons\:Ey\:Z  \advance\:y \:Y}
\def\:tDLn(#1,#2,#3,{\:tMvTo(#1,#2,#3)
                     \:tLnTo(}
\Define\:tMvTo(3){   \:X#1\:Xunits
   \:Y#2\:Yunits    \:Z#3\:Zunits}
\Define\:tMv(3){     \advance\:X  #1\:Xunits
   \advance\:Y  #2\:Yunits
   \advance\:Z  #3\:Zunits}
\Define\:tResize(3){   \:Xunits#1\:Xunits
   \:Yunits#2\:Yunits \:Zunits#3\:Zunits
   \:AdjRunits}
\Define\:tUnits(3){    \:Xunits#1  \relax
   \:Yunits#2  \relax       \:Zunits#3
   \:AdjRunits}\Define\:tLnTo(3){
   \:project   \edef\:temp{\:x\the\:x \:y\the\:y}
   \:X#1\:Xunits  \:Y#2\:Yunits  \:Z#3\:Zunits
   \:DLN}

\Define\:tLn(3){
   \:project   \edef\:temp{\:x\the\:x \:y\the\:y}
   \advance\:X  #1\:Xunits
   \advance\:Y  #2\:Yunits
   \advance\:Z  #3\:Zunits   \:DLN}

\Define\:DLN{
   \TwoDim      \:temp        \LineTo(\:x\du,\:y\du)
   \EndTwoDim }\Define\TwoDim{\bgroup\begingroup
   \def\EndTwoDim{\endgroup\egroup}
   \:loadIID
   \if:IIID  \:IIIDfalse \:project \:X\:x \:Y\:y
             \:CommonIID \fi
   \Units(\:Xunits,\:Yunits)}
\newif\if:rotated
\Define\RotatedAxes(2){\begingroup
   
   \if:IIID  \:IIIDfalse
             \:project \:X\:x \:Y\:y
             \:CommonIID    \fi
   \:DK#2\p@ \advance\:DK -#1\p@
\advance\:DK -\:CVXXX   \:NormalizeDeg\:DK
\ifdim \:DK<\:mmp \:wrn4{(#1,#2)}
\else
  \:DK#1\p@ \:NormalizeDeg\:DK    \:CosSin\:DK
  \:Xunitsx\:x  \:Xunitsy\:y
  \:DK#2\p@ \advance\:DK -\:XC \:NormalizeDeg\:DK \:CosSin\:DK
  \edef\Units(##1,##2){\noexpand\:Units(##1,##2)
    \:Xunitsx \:Cons\:Xunitsx\:Xunits
    \:Xunitsy \:Cons\:Xunitsy\:Xunits
    \:Yunitsx-\:Cons\:y\:Yunits
    \:Yunitsy \:Cons\:x\:Yunits  } \fi
 \def\MoveTo{\:rMvTo}
\def\Move{\:rMv}
\def\LineTo{\:rLnTo}
\def\Line{\:rLn}
\def\MoveF{\:rMvF}\def\:InitOval{\:rInitOval}

   \MarkLoc(:org) \Units(\:Xunits,\:Yunits)  \:rotatedtrue }
\newdimen\:Xunitsx  \newdimen\:Xunitsy
\newdimen\:Yunitsx  \newdimen\:Yunitsy
\Define\:rResize(2){
   \:Xunits #1\:Xunits    \:Yunits #2\:Yunits
   \:Xunitsx#1\:Xunitsx   \:Xunitsy#1\:Xunitsy
   \:Yunitsx#2\:Yunitsx   \:Yunitsy#2\:Yunitsy  }
\Define\:rMvTo{\MoveToLoc(:org) \:rMv}
\Define\:rMv(1){ \:rxy(#1)
   \advance\:X  \:x \advance\:Y  \:y}
\Define\:rLnTo(1){ \MarkLoc(:a) \:rMvTo(#1) {\LineToLoc(:a)} }
\Define\:rLn(1){ \:rxy(#1) \:Ln(\:x\du,\:y\du) }
\Define\:rxy(2){ \:Zunits\p@   \:Zunits#1\:Zunits
   \:x\:Cons\:Zunits\:Xunitsx
   \:y\:Cons\:Zunits\:Xunitsy
   \:Zunits\p@   \:Zunits#2\:Zunits
   \advance\:x  \:Cons\:Zunits\:Yunitsx
   \advance\:y  \:Cons\:Zunits\:Yunitsy}
\Define\:rMvF(1){  \:CosSin\:direction
   \:Zunits\p@            \:Zunits#1\:Zunits
   \:x\:Cons\:Zunits\:x   \:y\:Cons\:Zunits\:y
   \edef\:temp{(\:Cons\:x,\:Cons\:y)}   \expandafter\Move\:temp}
\Define\:AdjRunits{
   \:Xunitsx\:Xunits   \:Xunitsy\z@
   \:Yunitsx\z@        \:Yunitsy\:Yunits}
\Define\Text{  \setbox\:box
   \vtop\bgroup    \edef\DoReg{\the\DoReg}
      \hyphenpenalty\@M  \exhyphenpenalty\@M
      \catcode`\ 10 \catcode`\^^M13 \catcode`\^^I10
      \catcode`\&4  \let~\space
      \:Text}                       \catcode`\^^M13 %
\def\:Text(--#1--){%
      \:SetLines#1\hbox{}^^M--)^^M %
   \egroup                  %
   \if:IIID \TwoDim  \:DrawBox  \EndTwoDim %
   \else             \:DrawBox  \fi}       %
\def\:SetLines#1^^M{        %
   \def\:TextLine{#1}       %
   \ifx \:TextLine\:LastLine   \let\:temp\relax      %
   \else  \def\:temp{                                 %
             \:IndirectLines#1\relax~~--)~~\:SetLines}%
   \fi  \:temp }                      \catcode`\^^M9
\def\:IndirectLines#1~~{    \def\:TextLine{#1}
  \ifx \:TextLine\:LastLine   \let\:temp\relax
  \else  \def\:temp{\:AddLine{#1}\:IndirectLines}
  \fi  \:temp }

\Define\:LastLine{--)}

\def\:AddLine#1{
   \ifvmode \noindent    \hsize\z@ \else
      \hfil \penalty-500 \hbox{}    \fi
   \hfil#1
   \setbox\:box\hbox{#1}
   \ifdim \wd\:box>\hsize \hsize\wd\:box \fi}\def\TextPar#1#2{
   \def\:TxtPar##1(##2){##1(--##2--)}
   \edef\:temp{\expandafter\noexpand\csname :\string#2\endcsname}
   \edef#2{\noexpand\:TextPar\expandafter\noexpand\:temp}
   \expandafter\let\:temp\:undefined
   \expandafter#1\:temp}
                                   \catcode`\^^M13
\def\:TextPar#1{\begingroup     \catcode`\&4         %
   \catcode`\ 10 \catcode`\^^M13 \catcode`\^^I10 %
   \:TPar{#1}}                     \catcode`\^^M9  %

\def\:TPar#1(--#2--){\endgroup
   #1(--#2--)  }
 \newdimen\:direction  \newdimen\:Vdirection
\Define\:Rotate(1){   \advance \:direction   #1\p@
                      \:NormalizeDeg\:direction  }
\Define\:RotateTo(1){ \:direction  #1\p@
                      \:NormalizeDeg\:direction  }

\Define\:tRotate(2){
   \advance\:Vdirection  #2\p@
   \:NormalizeDeg\:Vdirection  \:Rotate(#1)}
\Define\:tRotateTo(2){
   \:Vdirection #2\p@
   \:NormalizeDeg\:Vdirection  \:RotateTo(#1)}
\Define\LineF(1){\MarkLoc(,)\MoveF(#1) {\LineToLoc(,)}}
\Define\:MvF(1){    \:CosSin\:direction
   \:d#1\:Xunits   \advance\:X  \:Cons\:d\:x
   \:d#1\:Yunits   \advance\:Y  \:Cons\:d\:y   }
\Define\MoveFToOval(2){
   \:NormalizeDeg\:direction \:CosSin\:direction  \:Zunits\p@
   \:d#2\:Zunits    \:x\:Cons\:d\:x
   \ifdim \:d=\z@  \:err5{...,\:Cons\:d} \fi
   \:d#1\:Zunits    \:y\:Cons\:d\:y
   \ifdim \:d=\z@  \:err5{\:Cons\:d,...} \fi     \:SearchDir
   \XSaveUnits
      \if:rotated
        \:Zunits\p@     \:Zunits#1\:Zunits
        \:x\:Cons\:Zunits\:Xunits
        \:Zunits\p@     \:Zunits#2\:Zunits
        \:y\:Cons\:Zunits\:Yunits
      \else
        \:x#1\:Xunits  \:y#2\:Yunits  \fi     \relax
      \Units(\:x,\:y)
      \:Vdirection\:direction  \:direction\:ddd
      \MoveF(\@ne)  \:direction\:Vdirection  \XRecallUnits  }
\Define\:tMvF(1){    \:CosSin\:Vdirection
   \:d #1\:Zunits   \advance\:Z  \:Cons\:y\:d
   \:xx#1\:Xunits     \:yy#1\:Yunits
   \:xx\:Cons\:x\:xx    \:yy\:Cons\:x\:yy
   \:CosSin\:direction
   \:x\:Cons\:xx\:x    \:y\:Cons\:yy\:y
   \advance\:X  \:x   \advance\:Y  \:y } \Define\CSeg{\:Opt[]\:CSeg1}
\def\:CSeg[#1]#2(#3,#4){   \MarkLoc($^)
   \MoveToLoc(#4) \:x\:X \:y\:Y
   \if:IIID  \:d\:Z  \fi   \MoveToLoc(#3)
   \advance\:x -\:X  \:x#1\:x
   \advance\:y -\:Y  \:y#1\:y
   \if:IIID   \advance\:d -\:Z  \:d#1\:d \fi
   \:t{#2}
   \edef\:temp{\the\:t(
      \expandafter\:Frac\the\:x\noexpand\:du,
      \expandafter\:Frac\the\:y\noexpand\:du \if:IIID ,
      \expandafter\:Frac\the\:d\noexpand\:du \fi)}
   \MoveToLoc($^)     \:temp}\Define\LSeg{\:Opt[]\:LSeg1}
\def\:LSeg[#1]#2(#3,#4){   \:distance(#3,#4)
   \:d#1\:d  \:t{#2}
   \edef\:temp{\the\:t(\expandafter\:Frac\the\:d\noexpand\:du)}
   \:temp}\Define\DSeg{\:Opt[]\:DSeg1}

\def\:DSeg[#1]#2(#3,#4){   \MarkLoc(^)
   \MoveToLoc(#4)  \:xxx\:X   \:yyy\:Y \:xxxx\:Z
   \MoveToLoc(#3)
   \advance\:xxx -\:X   \advance\:yyy -\:Y
   \ifdim \:AbsVal\:xxx<\:mmmp
      \ifdim \:AbsVal\:yyy<\:mmmp \:wrn5{#3,#4}
   \fi\fi
   \if:IIID
      \advance \:xxxx  -\:Z
      \:divide\:xxx\:Xunitsx
      \:divide\:yyy\:Yunitsy
      \:divide\:xxxx\:Zunits
      \:x\:xxx  \:y\:yyy   \:ddd\z@  \:length
      \:x\:d    \:y\:xxxx  \:SearchDir
      \:yyyy\:ddd
      \:x\:xxx  \:y\:yyy
   \else
      \:x  \:AbsVal\:Yunitsy
\:y  \:AbsVal\:Xunitsy  \ifdim \:y>\:x  \:x\:y  \fi
\:y  \:AbsVal\:Yunitsx  \ifdim \:y>\:x  \:x\:y  \fi
\:y  \:AbsVal\:Xunitsx  \ifdim \:y>\:x  \:x\:y  \fi
\:K  \:InCons\:x  \relax
      \ifnum \:K<\thr@@     \:K\@ne
\else \ifnum \:K<\sixt@@n   \:K4
\else \ifnum \:K<\:XC       \:K\sixt@@n
\else \ifnum \:K<\@m        \:K\@cclvi
\fi \fi \fi \fi
\divide\:xxx \:K
\divide\:yyy \:K
      \:x \:Cons\:Yunitsy\:xxx
\advance\:x -\:Cons\:Yunitsx\:yyy
      \:y-\:Cons\:Xunitsy\:xxx
\advance\:y \:Cons\:Xunitsx\:yyy
   \fi
   \:SearchDir   \:ddd#1\:ddd   \:t{#2}
   \edef\:temp{\the\:t(
      \:Cons\:ddd \if:IIID ,\:Cons\:yyyy \fi)}
   \MoveToLoc(^) \:temp}   \def\:theDoReg{\def\DoReg{\the\:DoReg}}

\Define\MarkLoc{  \:theDoReg
   \expandafter\edef \csname \:MarkLoc}

\Define\MarkGLoc{  \:theDoReg
   \expandafter\xdef \csname \:MarkLoc}

\Define\:MarkLoc(1){ Loc\space#1:\endcsname{
   \:X\the\:X  \:Y\the\:Y
   \if:IIID \:Z\the\:Z \fi  }       \let\DoReg\:DoReg }

\Define\MoveToLoc(1){   \:theDoReg   \expandafter\ifx
   \csname Loc\space#1:\endcsname\relax \:err2{#1}\fi
   \csname Loc\space#1:\endcsname     \let\DoReg\:DoReg  }

\Define\MarkPLoc(1){  \:theDoReg
   \if:IIID   \:project
      \expandafter\edef \csname Loc\space#1:\endcsname{
         \:X\the\:x  \:Y\the\:y  \:Z\z@}
   \else \:err1\MarkPLoc  \fi   \let\DoReg\:DoReg}

\Define\WriteLoc(1){{ \:theDoReg \edef\:temp{#1}
   \ifx \:temp\empty \else \MoveToLoc(#1) \fi
   \immediate\write\sixt@@n{...
      \:temp=(\the\:X,\the\:Y\if:IIID,\the\:Z\fi)}}}

\Define\:LnToLoc(1){
   \:x\:X \:y\:Y   \MoveToLoc(#1)
   { \:LnTo(\:x\du,\:y\du) }}

\Define\:tLnToLoc(1){
   \:xx\:X \:yy\:Y \:dd\:Z     \MoveToLoc(#1)
   { \:tLnTo(\:xx\du,\:yy\du,\:dd\du) }}\def\:GetLine#1#2#3#4#5{
   \MoveToLoc(#1)   \divide\:X \:eight  \divide\:Y \:eight
   #3\:X  #4\:Y
   \MoveToLoc(#2)   \divide\:X \:eight  \divide\:Y \:eight
   \advance #3 -\:X  \advance #4 -\:Y
   #5\:Cons#3\:Y      \advance #5 -\:Cons#4\:X
   \divide #3 \:eight     \divide  #4 \:eight \relax      }
\def\MoveToLL(#1,#2)(#3,#4){
   \:GetLine{#1}{#2}\:x \:y \:xxx
   \:GetLine{#3}{#4}\:xx\:yy\:xxxx
   \:ddd \:Cons\:x \:yy     \advance\:ddd -\:Cons\:xx\:y
   \ifdim  \:AbsVal\:ddd < \:mmmp
      \:X\@cclv\p@  \:Y\:X
      \:wrn3{(\string#1,\string#2)(\string#3,\string#4)}
   \else
      \:divide\:xxx\:ddd    \:divide\:xxxx\:ddd
      \:X\:Cons\:xxx\:xx   \advance\:X -\:Cons\:xxxx\:x
      \:Y\:Cons\:xxx\:yy   \advance\:Y -\:Cons\:xxxx\:y
   \fi   }\Define\MoveToCC{\:Opt[]\:MoveToCC{}}

\def\:MoveToCC[#1](#2,#3)(#4,#5){
  \:UserUnits(#2,#3)(#4,#5)
  \:distance($#2,$#4)
\ifnum \:d<\:mp \MoveToLoc(#3)
   \:wrn3{(#1,#2)(#3,#4)}
\else                 \:xx \:d
   \:distance($#2,$#3)  \:xxx\:d
   \:distance($#5,$#4)
     \:yy \:xxx  \advance\:yy -\:d
\:yyy\:xxx  \advance\:yyy \:d
\:divide\:yy\:xx     \:yy\:Cons\:yy\:yyy
\advance\:yy \:xx  \divide\:yy \tw@
     \:yyy \ifdim \:AbsVal\:xxx>\:AbsVal\:yy \:xxx \else \:yy \fi
\ifdim \:AbsVal\:yyy<\:mp \:yyy\z@ \else
   \:divide\:xxx\:yyy  \:sqr\:xxx
   \:yyyy\:yy  \:divide\:yyyy\:yyy  \:sqr\:yyyy
   \advance\:xxx -\:yyyy   \:Sqrt\:xxx
   \:yyy\:Cons\:yyy\:xxx
\fi
     \MoveToLoc($#4)  \:x\:X  \:y\:Y  \:divide\:yy\:xx
\MoveToLoc($#2)
\advance\:x -\:X     \advance\:y -\:Y
\advance\:X \:Cons\:yy\:x
\advance\:Y \:Cons\:yy\:y
     \:divide\:yyy\:xx
\advance\:X  #1\:Cons\:yyy\:y
\advance\:Y -#1\:Cons\:yyy\:x

  \fi  \:SysUnits  }\def\:UserUnits(#1,#2)(#3,#4){
   \:xx\:Xunitsx  \:xxx\:Xunitsy
   \:yy\:Yunitsx  \:yyy\:Yunitsy
   \:xxxx\:Cons\:yy\:xxx  \advance\:xxxx \:Cons\:yyy\:xx
   \ifdim \:AbsVal\:xxxx>\:mmmp
      \:divide\:xx\:xxxx   \:divide\:xxx{-\:xxxx}
      \:divide\:yy{\:xxxx} \:divide\:yyy\:xxxx
   \fi
   \:UnLoc(#1)  \:UnLoc(#2)
   \:UnLoc(#3)  \:UnLoc(#4)}\Define\:SysUnits{
   \MoveTo(\:Cons\:X,\:Cons\:Y)}\Define\:UnLoc(1){
   \MoveToLoc(#1)  \:d\:Cons\:yyy\:X
   \advance\:d \:Cons\:yy\:Y
   \:Y\:Cons\:xx\:Y
   \advance\:Y  \:Cons\:xxx\:X   \:X\:d
   \MarkLoc($#1)}

\def\:MoveToLC[#1](#2,#3)(#4,#5){
   \:UserUnits(#2,#3)(#4,#5)
   \MoveToLoc($#2)  \:x\:X  \:y\:Y
\MoveToLoc($#3)  \advance\:x -\:X
                \advance\:y -\:Y
   \edef\:temp{ \:xxx\the\:x  \:yyy\the\:y }
\MoveToLoc($#4)
\advance\:X \:y  \advance\:Y -\:x
\MarkLoc(^$) \MoveToLL($#4,^$)($#2,$#3)
\MarkLoc(^$)
   \:distance($#4,$#5)  \:xx\:d  \:distance($#4,^$)
\ifdim      \:d>\:xx        \:wrn3{(#2,#3)(#4,#5)}
\else \ifdim \:d<\:mmp     \:yy\:xx   \else   \:yy\:d
      \:divide\:yy\:xx  \:sqr\:yy
      \:yy-\:yy   \advance\:yy \p@
      \:Sqrt\:yy  \:yy\:Cons\:xx\:yy
\fi \fi
   \:temp   \:x\:xxx  \:y\:yyy  \:length  \:xx\:d
\:divide\:xxx\:xx
\:divide\:yyy\:xx
\advance\:X  #1\:Cons\:yy\:xxx
\advance\:Y  #1\:Cons\:yy\:yyy
   \:SysUnits }  \def\Object#1{\:Opt(){\:DefineSD#1}0}

\def\:DefineSD#1(#2){\begingroup  \:multid#1
   \:DraCatCodes   \:DefSD#1(#2)}

\def\:DefSD#1(#2)#3{
   \expandafter\::Define\csname\string#1.\endcsname(#2){
      \:t{\:SubD{#3}}
      \if:IIID \edef\:temp{\noexpand\TwoDim \the\:t
                           \noexpand\EndTwoDim}
      \else    \def\:temp{\the\:t}   \fi         \:temp}
   \def#1{\def\:SDname{\csname\string#1.\endcsname}
          \:Opt[]\:CallSD{}}}

\def\:CallSD[#1]{ \edef\:Entry{#1} \:SDname }\def\:SubD#1{
   \let\::RecallXLoc\:AddXLoc   \gdef\:AddXLoc{}
   \edef\:RecallBor{ \global\:LBorder  \the\:LBorder
                  \global\:RBorder  \the\:RBorder
                  \global\:TeXLoc\the\:TeXLoc }
\global\:TeXLoc\z@
\setbox\:box\vbox{\EntryExit(0,0,0,0)
\begingroup
   
   \:InDraw  #1
\endgroup
\:SetDrawWidth                    \let\:XLoc\relax
\xdef\:AddXLoc{\:dd\the\:LBorder  \:AddXLoc}}
\:RecallBor
   \:ddd\dp\:box
   \ifx \:Entry\empty
   \:DrawBox
\else
   \let\:RecallIn\:InOut
   \:x\:X    \:y\:Y
   \def\:XLoc(##1,##2,##3){
      \def\:temp{##1}
      \ifx \:temp\:Entry \:X\:x  \advance\:X -##2
                         \:Y\:y  \advance\:Y -##3
      \fi}
   \:AddXLoc
   \advance\:X  \:dd      \advance\:Y -\:ddd
   \EntryExit(-1,-1,\:InOut2,\:InOut3)  \:DrawBox
   \let\:InOut\:RecallIn
\fi
   \MarkLoc(^)
      \MoveToExit(-1,-1)
   \:xxx\:X   \:yyy\:Y   \advance\:yyy  \:ddd
   \def\:XLoc(##1,##2,##3){
      \:X\:xxx  \advance\:X  ##2     \advance\:X -\:dd
      \:Y\:yyy  \advance\:Y  ##3
      \MarkLoc(##1)}
   \:AddXLoc
   \MoveToLoc(^)
   \ifx \:Entry\empty \else     \MoveToLoc(\:Entry) \fi
   \global\let\:AddXLoc\::RecallXLoc }
\Define\:MarkXLoc(1){  \:theDoReg
   \let\:XLoc\relax
   \xdef\:AddXLoc{\:AddXLoc \:XLoc(#1,\the\:X,\the\:Y)}
   \let\DoReg\:DoReg}
   \catcode`\ 10 \catcode`\^^M5 \catcode`\^^I10
   \let\wlog\:wlog  \let\:wlog\:undefined
\:RestoreCatcodes     \tracingstats1
\newtheorem{thm}{Theorem}[section]
\newtheorem{exa}[thm]{Example}

\newtheorem{defn}[thm]{Definition}
\newtheorem{prop}[thm]{Proposition}
\newtheorem{rem}[thm]{Remark}

\newcommand{\B}{\mathbb{B}}
\newcommand{\RR}{\mathbb{R}}
\newcommand{\E}{\mathbb{S}}
\newcommand{\Z}{\mathbb{Z}}
\title{Blanchfield and Seifert algebra in high-dimensional knot theory}
\author{Andrew Ranicki}
\address{School of Mathematics\newline
\indent University of Edinburgh\newline
\indent King's Buildings\newline
\indent Edinburgh EH9 3JZ \newline
\indent Scotland, UK\hfill e-mail:{\tt aar@maths.ed.ac.uk}}

\begin{document}

\maketitle

\hfill{\it Dedicated to S.P.Novikov}

\bigskip

Novikov \cite{Nov} initiated the study of the algebraic properties
of quadratic forms over polynomial extensions by a far-reaching
analogue of the Pontrjagin-Thom transversality construction of a
Seifert surface of a knot and the infinite cyclic cover of the
knot exterior.  In this paper the analogy is applied to explain
the relationship between the Seifert forms over a ring with
involution $A$ and Blanchfield forms over the Laurent polynomial
extension $A[z,z^{-1}]$.

The rings $A$ and $A[z,z^{-1}]$ correspond to the two ways of
associating algebraic invariants to an $n$-knot $k:S^n \subset
S^{n+2}$ with $A=\Z$~:
\begin{itemize}
\item[(i)] The $\Z[z,z^{-1}]$-module invariants of the canonical infinite
cyclic cover $\overline{M}=p^*\RR$ of the exterior of $k$
$$M^{n+2}~=~{\rm cl.}(S^{n+2}\backslash k(S^n)\times D^2) \subset S^{n+2}$$
with $k(S^n)\times D^2 \subset S^{n+2}$ a regular
neighbourhood of $k(S^n)$ in $S^{n+2}$, $p:M \to S^1$ a map
inducing an isomorphism $p^*:H^1(S^1) \cong H^1(M)$, and
$\partial M = S^n \times S^1$.
\item[(ii)] The $\Z$-module invariants of a codimension 1 submanifold
$N^{n+1} \subset S^{n+2}$ with boundary
$$\partial N~=~k(S^n) \subset S^{n+2}~,$$
i.e. a Seifert surface for $k$.
\end{itemize}
The knot $k$ has a unique exterior $M$, and many Seifert surfaces $N$.
For any $p:M \to S^1$ which is transverse regular at $1 \in S^1$
the inverse image
$$N~=~p^{-1}(1) \subset M$$
is a Seifert surface for $k$.
Conversely, any $N$ can be used to construct $\overline{M}$ as an
infinite union of fundamental domains $(M_N;N,zN)$
$$\overline{M}~=~\bigcup\limits^{\infty}_{j=-\infty}z^jM_N~.$$

Chapter \ref{modules} deals with the following concepts~:
\begin{itemize}
\item[(i)] A {\it Seifert module} over $A$ is a pair
$$(P,e)~=~(~\hbox{f.g. projective $A$-module}~,~\hbox{endomorphism}~)~.$$
\item[(ii)] A {\it Blanchfield module} $B$ is a homological
dimension 1  $A[z,z^{-1}]$-module such that $1-z:B \to B$ is an automorphism.
\item[(iii)] The {\it covering} of a Seifert module $(P,e)$ is the
Blanchfield module
$$B(P,e)~=~{\rm coker}(1-e+ze:P[z,z^{-1}] \to P[z,z^{-1}])~.$$
\end{itemize}
The covering construction $B:(P,e) \mapsto B(P,e)$ is an algebraic
version of the construction of the infinite cyclic cover $\overline{M}$
from $(M_N;N,zN)$.  Theorem \ref{cover1} proves that every Blanchfield
module $B$ is isomorphic to the covering $B(P,e)$ of a Seifert module
$(P,e)$.  Moreover, morphisms of Blanchfield modules are characterized
in terms of morphisms of Seifert modules.

Chapter \ref{kernel} characterizes the Seifert modules $(P,e)$ such
that $B(P,e)=0$, and also the morphisms of Seifert
modules with covering an isomorphism of Blanchfield modules.

Chapter \ref{forms} deals with the following concepts, where $\eta=\pm
1$, and $A$ is a ring with involution~:
\begin{itemize}
\item[(i)] An {\it $\eta$-symmetric Seifert form} $(P,\theta)$ is a f.g. projective
$A$-module $P$ together with an $A$-module morphism
    $$\theta~:~P \to P^*~=~{\rm Hom}_A(P,A)$$
such that $\theta+\eta\theta^*:P \to P^*$ is an isomorphism.
\item[(ii)] An {\it $\eta$-symmetric Blanchfield form} $(B,\phi)$ is a Blanchfield
$A[z,z^{-1}]$-module $B$ together with an isomorphism
$$\phi~:~B \to B\widehat{~}~=~{\rm Ext}^1_{A[z,z^{-1}]}(B,A[z,z^{-1}])$$
such that $\eta\phi\widehat{~}=\phi$.
\item[(iii)] The {\it covering} of a $(-\eta)$-symmetric Seifert form $(P,\theta)$ is the
$\eta$-symmetric Blanchfield form
$$B(P,\theta)~=~(B(P,e),\phi)$$
with $e=(\theta-\eta\theta^*)^{-1}\theta:P \to P$ and
$\phi=(1-z^{-1})\zeta_{(P,e)}B(\theta-\eta\theta^*)$ (see \ref{Bcover}
for details).
\end{itemize} Theorem \ref{cover2} gives
an algorithmic proof that every $\eta$-symmetric Blanchfield form
$(B,\phi)$ over $A[z,z^{-1}]$ is isomorphic to the covering
$B(P,\theta)$ of a $(-\eta)$-symmetric Seifert form $(P,\theta)$ over $A$.

Chapter \ref{Witt} deals with algebraic $L$-theory.  Theorem
\ref{cover3} identifies the Witt group of $\eta$-symmetric
Blanchfield forms over $A[z,z^{-1}]$ with the Witt group of
$(-\eta)$-symmetric Seifert forms over $A$.  Theorem \ref{invert}
identifies this group with a quotient of the Witt group of
$\eta$-symmetric forms over the universal localization
$\Pi^{-1}A[z,z^{-1},(1-z)^{-1}]$ of $A[z,z^{-1}]$ inverting $1-z$
and the set $\Pi$ of $A$-invertible matrices over $A[z,z^{-1}]$.
For $A=\Z$, $\eta=(-1)^{i+1}$, $i \geqslant 2$ this is an expression of the
$(2i-1)$-dimensional knot cobordism group as
$$C_{2i-1}~=~{\rm coker}(L_{2i+2}(\Z[z,z^{-1},(1-z)^{-1}]) \to
L_{2i+2}(P^{-1}\Z[z,z^{-1},(1-z)^{-1}]))$$
with $P=\{p(z)\vert p(1)=1\} \subset \Z[z,z^{-1}]$ the Alexander polynomials.

I am grateful to the Mathematics Department of the University of California,
San Diego, which I was visiting January--March 2001 when work on
this paper was started.  I am also grateful to Peter Teichner and Des
Sheiham for various conversations and e-mails. In particular, Des simplified
the formulation of Proposition 3.9 (iii).

\section{Blanchfield and Seifert modules}\label{modules}

Let $A$ be a ring, with Laurent polynomial extension $A[z,z^{-1}]$.

\begin{defn} {\rm
A f.g. projective $A[z,z^{-1}]$-module is {\it induced}
if it is of the form
$$P[z,z^{-1}]~=~A[z,z^{-1}]\otimes_AP$$
for a f.g. projective $A$-module $P$.\hfill\qed}
\end{defn}

We shall make frequent use of the identity
$${\rm Hom}_{A[z,z^{-1}]}(P[z,z^{-1}],Q[z,z^{-1}])~=~
{\rm Hom}_A(P,Q)[z,z^{-1}]$$
with $P,Q$ f.g. projective $A$-modules.

\begin{defn} \label{Bdefn}
{\rm (i) A {\it Blanchfield $A[z,z^{-1}]$-module} $B$ is an $A[z,z^{-1}]$-module
such that
\begin{itemize}
\item[(a)] $1-z:B \to B$ is an automorphism,
\item[(b)] there exists an induced f.g. projective $A[z,z^{-1}]$-module resolution
$$0 \to P_1[z,z^{-1}] \xymatrix{\ar[r]^d&} P_0[z,z^{-1}] \to B \to 0~.$$
\end{itemize}
(ii) The {\it Blanchfield category} $\B(A[z,z^{-1}])$ has objects
Blanchfield $A[z,z^{-1}]$-modules and $A[z,z^{-1}]$-module morphisms.}
\hfill\qed \end{defn}

Write the augmentation as
$$\epsilon~:~A[z,z^{-1}] \to A~;~z \mapsto 1~.$$

\begin{prop} \label{Bprop}
Let $C$ be a 1-dimensional induced f.g. projective $A[z,z^{-1}]$-module
chain complex with
$$d~=~\sum\limits^k_{j=0}d_jz^j~:~C_1~=~P_1[z,z^{-1}] \to C_0~=~P_0[z,z^{-1}]~.$$
The homology $A[z,z^{-1}]$-module
$$B~=~H_0(C)~=~{\rm coker}(d)$$
is a Blanchfield module if and only if the $A$-module morphism
$$\epsilon(d)~=~\sum\limits^k_{j=0}d_j~:~P_1 \to P_0$$
is an isomorphism.
\end{prop}
\begin{proof}
If $B$ is a Blanchfield module the inverse isomorphism
$(1-z)^{-1}:B \to B$ is resolved by an $A[z,z^{-1}]$-module
chain map $f:C \to C$
$$\xymatrix{
0 \ar[r] &C_1 \ar[r]^d \ar[d]^{f_1} & C_0 \ar[d]^{f_0} \ar[r]
& B \ar[r] \ar[d]^{(1-z)^{-1}}  & 0 \\
0 \ar[r] & C_1 \ar[r]^{d} & C_0 \ar[r] & B \ar[r] &0}$$
so that $f:C \to C$ is chain homotopy inverse to $1-z:C \to C$.
A chain homotopy
$$g~:~f(1-z)~ \simeq~ 1~:~C \to C$$
is defined by an $A[z,z^{-1}]$-module morphism
$$g~=~\sum\limits^s_{j=r}z^jg_j~:~C_0~=~P_0[z,z^{-1}] \to C_1~=~P_1[z,z^{-1}]$$
such that
$$\begin{array}{l}
1- f_0(1-z) ~=~dg ~:~C_0~=~P_0[z,z^{-1}] \to C_0~=~P_0[z,z^{-1}]~,\\[1ex]
1- f_1(1-z) ~=~gd ~:~C_1~=~P_1[z,z^{-1}] \to C_1~=~P_1[z,z^{-1}]~,
\end{array}$$
and
$$\epsilon(g)~=~\sum\limits^s_{j=r}g_j~:~P_0 \to P_1$$
is an $A$-module isomorphism inverse to $\epsilon(d)=\sum\limits^k_{j=0}d_j:P_1 \to P_0$.

Conversely, suppose that $\epsilon(d):P_1 \to P_0$ is an
isomorphism, with inverse
$$\epsilon(d)^{-1}~=~h~:~P_0 \to P_1~,$$
so that
$$\begin{array}{l}
1-dh ~=~\sum\limits^k_{j=0}(1-z^j)d_jh ~:~C_0~=~P_0[z,z^{-1}] \to C_0~=~P_0[z,z^{-1}]~,\\[1ex]
1-hd ~=~\sum\limits^k_{j=0}(1-z^j)hd_j ~:~C_1~=~P_1[z,z^{-1}] \to C_1~=~P_1[z,z^{-1}]~.
\end{array}$$
The $A[z,z^{-1}]$-module morphisms
$$\begin{array}{l}
f_0 ~=~(1-dh)(1-z)^{-1} ~:~C_0~=~P_0[z,z^{-1}] \to C_0~=~P_0[z,z^{-1}]~,\\[1ex]
f_1 ~=~(1-hd)(1-z)^{-1} ~:~C_1~=~P_1[z,z^{-1}] \to C_1~=~P_1[z,z^{-1}]~.
\end{array}$$
are the components of a chain equivalence $f:C \to C$ chain homotopy
inverse to $1-z:C \to C$, with a chain homotopy
$$h~:~f(1-z)~ \simeq~ 1~:~C \to C~.$$
It remains to verify that $d:C_1 \to C_0$ is injective.
If $x \in {\rm ker}(d:C_1 \to C_0)$ then
$$x~=~(1-hd)(x)~=~(1-z)f_1(x) \in C_1$$
with $f_1(x) \in {\rm ker}(d:C_1 \to C_0)$ by the injectivity
of $1-z:C_0 \to C_0$. It follows that for any integer $j \geqslant 1$
$$x~=~(1-z)^j(f_1)^j(x) \in C_1~,$$
and
$$x \in \bigcap^{\infty}_{j=1}(1-z)^j(C_1)~=~\{0\}\subset C_1~.$$
\end{proof}

Proposition \ref{Bprop} is the special case $n=0$ of~:

\begin{prop} {\it The following conditions on an $(n+1)$-dimensional
induced f.g.  projective $A[z,z^{-1}]$-module chain complex $C$ are
equivalent~:
\begin{itemize}
\item[\rm (i)] there exists a homology equivalence $C \to B$ to an
$n$-dimensional chain complex in the Blanchfield category
$\B(A[z,z^{-1}])$,
\item[\rm (ii)] $H_*(A\otimes_{A[z,z^{-1}]}C)=0$.
\end{itemize}}
\end{prop}
\begin{proof} As for Proposition 3.1.2 of Ranicki \cite{RESATS}.
\end{proof}

\begin{exa} \label{Bexa}
 {\rm Let $M$ be a finite $CW$ complex with a homology equivalence $p:M
\to S^1$, such as a knot complement. Let $\overline{M}=p^*\RR$ be
the pullback infinite cyclic cover of $M$, and let
$C=C(\overline{p}:\overline{M} \to \RR)_{*+1}$ be the relative
cellular $\Z[z,z^{-1}]$-module chain complex of the induced
$\Z$-equivariant cellular map $\overline{p}:\overline{M} \to \RR$,
with $H_*(C)=\widetilde{H}_*(\overline{M})$ the reduced homology
of $\overline{M}$. Then $H_*(\Z\otimes_{\Z[z,z^{-1}]}C)=H_{*+1}
(\overline{p})=0$, and $C$ is homology equivalent to a finite
chain complex in the Blanchfield category
$\B(\Z[z,z^{-1}])$.\hfill$\qed$}
\end{exa}

\begin{defn} \label{Sdefn}
{\rm (i) A {\it Seifert $A$-module} $(P,e)$ is a f.g. projective $A$-module
$P$ together with an endomorphism $e:P \to P$.\\
(ii) A {\it morphism} of Seifert $A$-modules
$g:(P,e) \to (P',e')$
is an $A$-module morphism such that
$$e'g~=~ge ~:~P \to P'$$
(iii) The {\it Seifert category} $\E(A)$ has objects
Seifert $A$-modules and morphisms as in (ii). \hfill\qed}
\end{defn}

Seifert modules determine Blanchfield modules by :

\begin{defn} {\rm (i) The {\it covering} of a Seifert $A$-module $(P,e)$ is
the Blanchfield $A[z,z^{-1}]$-module
$$B(P,e)~=~{\rm coker}(1-e+ze:P[z,z^{-1}] \to P[z,z^{-1}])$$
with the resolution
$$C(P,e)~:~C_1~=~P[z,z^{-1}] \xymatrix{\ar[r]^{1-e+ze}&}
C_0~=~P[z,z^{-1}]~.$$
(ii) The {\it covering} of a Seifert $A$-module morphism
$g:(P,e) \to (P',e')$ is
the Blanchfield $A[z,z^{-1}]$-module morphism
$$B(g)~:~B(P,e) \to B(P',e')~;~x \mapsto g(x)$$
resolved by the chain map
$$\xymatrix@R+2pt{~C(P,e): \ar[d]_{C(g)} \\~C(P',e'): }
\xymatrix@C+25pt{P[z,z^{-1}] \ar[r]^{1-e+ze}
\ar[d]^g & P[z,z^{-1}] \ar[d]^g \\
P'[z,z^{-1}] \ar[r]^{1-e'+ze'} & P'[z,z^{-1}]~.}$$
\hfill\qed}
\end{defn}

\begin{thm} \label{cover1}
The covering construction defines a functor of additive categories
$$B~:~\E(A) \to \B(A[z,z^{-1}])~;~(P,e) \mapsto B(P,e)$$
such that
\begin{itemize}
\item[(i)] Every Blanchfield $A[z,z^{-1}]$-module $B$ is isomorphic to
the covering $B(P,e)$ of a Seifert $A$-module $(P,e)$.
\item[(ii)] The coverings of $e,1-e:(P,e) \to (P,e)$ are automorphisms
$$B(e)~=~(1-z)^{-1}~~,~~B(1-e)~=~(1-z^{-1})^{-1}~:~B(P,e) \to B(P,e)~,$$
with inverses
$$B(e)^{-1}~=~1-z~~,~~B(1-e)^{-1}~=~1-z^{-1}~:~B(P,e) \to B(P,e)~.$$
\item[(iii)] Every morphism of Blanchfield $A[z,z^{-1}]$-modules
$f:B(P,e)\to B(P',e')$ is of the type
$$f~=~B(g)t^{-k}$$
for some morphism of Seifert $A$-modules $g:(P,e) \to (P',e')$ and $k \geqslant 0$,
with $t$ the automorphism
$$t~=~B(e(1-e))~=~((1-z)(1-z^{-1}))^{-1}~:~B(P,e) \to B(P,e)~.$$
\item[(iv)] Two morphisms $g_1,g_2:(P,e) \to (P',e')$ are such that
$$B(g_1)t^{-k_1}~=~B(g_2)t^{-k_2}~:~B(P,e) \to B(P',e')$$
for some $k_1,k_2 \geqslant 0$ if and only if
$$(g_1(e(1-e))^{k_2}-g_2(e(1-e))^{k_1})(e(1-e))^{\ell}~=~0~:~P \to P'$$
for some $\ell \geqslant 0$.
\end{itemize}
\end{thm}
\begin{proof}
(i) By Proposition \ref{Bprop} it may be assumed that $B=H_0(C)$ with
$$d~=~\sum\limits^k_{j=0}d_jz^j~:~C_1~=~P_1[z,z^{-1}] \to C_0~=~P_0[z,z^{-1}]$$
for f.g. projective $A$-modules $P_0,P_1$, such that the augmentation
$A$-module morphism
$$\epsilon(d)~=~\sum\limits^k_{j=0}d_j~:~P_1 \to P_0$$
is an $A$-module isomorphism.

Let $s$ be another indeterminate over $A$, and use the isomorphism of rings
$$A[s,s^{-1},(1-s)^{-1}] \to A[z,z^{-1},(1-z)^{-1}]~;~s \mapsto (1-z)^{-1}$$
to identify
$$A[s,s^{-1},(1-s)^{-1}]~=~A[z,z^{-1},(1-z)^{-1}]~,$$
with
$$s~=~(1-z)^{-1}~~,~~z~=~s^{-1}(s-1)~.$$
The $A[z,z^{-1},(1-z)^{-1}]$-module morphism induced by
$d:C_1 \to C_0$
$$d~=~\sum\limits^k_{j=0}d_jz^j~:~
P_1[z,z^{-1},(1-z)^{-1}] \to P_0[z,z^{-1},(1-z)^{-1}]$$
is expressed in terms of $s$ as
$$d~=~\sum\limits^k_{j=0}(s^{-1}(s-1))^jd_j~:~
P_1[s,s^{-1},(1-s)^{-1}] \to P_0[s,s^{-1},(1-s)^{-1}]~.$$
The $A[s]$-module morphism
$$\Delta~=~
\sum\limits^k_{j=0}d_j\epsilon(d)^{-1}s^{k-j}(s-1)^j~:~P_0[s] \to P_0[s]$$
induces the $A[s,s^{-1},(1-s)^{-1}]$-module morphism
$$\Delta~=~s^kd\epsilon(d)^{-1}~:~
P_0[s,s^{-1},(1-s)^{-1}] \to P_0[s,s^{-1},(1-s)^{-1}]~.$$
Now $\Delta$ is a degree $k$ polynomial with coefficients
$\Delta_j \in {\rm Hom}_A(P_0,P_0)$
$$\Delta~=~\sum\limits^k_{j=0}\Delta_js^j~:~P_0[s] \to P_0[s]$$
such that $\Delta_k=1$. The Seifert $A$-module
$$(P,e)~=~\big(P_0 \oplus P_0 \oplus \dots \oplus P_0~(k~{\rm terms}),
\begin{pmatrix}
0 & 0 & 0 & \dots & -\Delta_0 \cr
1 & 0 & 0 & \dots & -\Delta_1 \cr
0 & 1 & 0 & \dots & -\Delta_2 \cr
\vdots & \vdots &\vdots &\ddots & \vdots \cr
0 & 0 & 0 & \dots & -\Delta_{k-1}
\end{pmatrix}\big)$$
is such that there is defined an exact sequence of $A[s]$-modules
$$0\to P_0[s] \xymatrix{\ar[r]^{\Delta} &} P_0[s] \to  P \to 0$$
with $s$ acting on $P$ by $e$ and
$$P_0[s] \to P~:~\sum\limits^{\infty}_{j=0}s^jx_j \mapsto
(x_0,x_1,\dots,x_{k-1})~.$$
The covering of $(P,e)$ is the induced $A[s,s^{-1},(1-s)^{-1}]$-module
$$B(P,e)~=~A[s,s^{-1},(1-s)^{-1}]\otimes_{A[s]}P$$
and the isomorphism of exact sequences of $A[s,s^{-1},(1-s)^{-1}]$-modules
$$\xymatrix@C+5pt{0 \ar[r] & P_1[s,s^{-1},(1-s)^{-1}] \ar[r]^d
\ar[d]_{\epsilon(d)}
&P_0[s,s^{-1},(1-s)^{-1}] \ar[r] \ar[d]_{s^k}  & B \ar[d] \ar[r] & 0 \\
0 \ar[r] & P_0[s,s^{-1},(1-s)^{-1}] \ar[r]^{\Delta} & P_0[s,s^{-1},(1-s)^{-1}]
\ar[r] & B(P,e) \ar[r] & 0}$$
includes an isomorphism
$$B~ \cong~ B(P,e)~.$$
(ii)  The $A[z,z^{-1}]$-module chain maps
$$C(e)~,~1-z~:~C(P,e) \to C(P,e)$$
are inverse chain homotopy equivalences, with
$$(1-z)C(e)~=~C(e)(1-z)~:~C(P,e) \to C(P,e)$$
and a chain homotopy
$$1~:~(1-z)C(e)~ \simeq~ {\rm id}~:~C(P,e) \to C(P,e)~.$$
Likewise, the $A[z,z^{-1}]$-module chain maps
$$C(1-e)~,~-z^{-1}(1-z)^{-1}~:~C(P,e) \to C(P,e)$$
are inverse chain homotopy equivalences, with
$$-z^{-1}(1-z)C(1-e)~=~C(1-e)(-z^{-1}(1-z))~:~C(P,e) \to C(P,e)$$
and a chain homotopy
$$z^{-1}~:~-z^{-1}(1-z)C(1-e)~ \simeq~ {\rm id}~:~C(P,e) \to C(P,e)~.$$
(iii) With $s=(1-z)^{-1}$ as in (i) define
$$t~=~s(1-s)~=~-z(1-z)^{-2}~,$$
and identify
$$A[s,s^{-1},(1-s)^{-1}]~=~A[s,t^{-1}]~=~A[z,z^{-1},(1-z)^{-1}]~,$$
Suppose given Seifert $A$-modules $(P,e)$, $(P',e')$ and
a morphism of Blanchfield $A[z,z^{-1}]$-modules
$f:B(P,e) \to B(P',e')$. Resolve $f$ by an $A[s,t^{-1}]$-module chain map
$$\xymatrix@C+25pt{C_1=P[s,t^{-1}] \ar[r]^{s-e}
\ar[d]^{f_1} & C_0=P[s,t^{-1}] \ar[d]^{f_0} \\
C'_1=P'[s,t^{-1}] \ar[r]^{s-e'} & C'_0=P'[s,t^{-1}]}$$
with
$$\begin{array}{l}
f_0~=~t^{-k}\sum\limits^{\ell}_{j=0}s^jf_{0,j}~:~P[s,t^{-1}] \to P'[s,t^{-1}]~,\\[1ex]
f_1~=~t^{-k}\sum\limits^{\ell}_{j=0}s^jf_{1,j}~:~P[s,t^{-1}] \to P'[s,t^{-1}]
\end{array}$$
for some $A$-module morphisms $f_{0,j},f_{1,j}:P \to P'$.
The morphism of Seifert $A$-modules
$$g~:~(P,e) \to (P',e')$$
with
$$g~=~\sum\limits^k_{j=0}(e')^jf_{0,j}~:~P \to P'$$
is such that
$$f~=~B(g)t^{-k}~:~B(P,e) \to B(P',e')$$
with
$$t~=~B(e(1-e))~:~B(P,e) \to B(P,e)~.$$
(Example: $-z=(1-e)^2t^{-1}:B(P,e)\to B(P,e)$.)\\
(iv) It suffices to show that a morphism of Seifert $A$-modules $g:(P,e) \to (P',e')$
is such that
$$B(g)~=~0~:~B(P,e) \to B(P',e')$$
if and only if for some $k \geqslant 0$
$$g(e(1-e))^k~=~0~:~P \to P~.$$
Now $B(g)=0$ if and only if there exists an $A[z,z^{-1}]$-module chain homotopy
$$h~:~g~ \simeq~ 0~:~C(P,e) \to C(P',e')$$
with
$$\xymatrix@R+2pt{~C(P,e): \ar[d]_g \\~C(P',e'): }
\xymatrix@C+25pt{P[z,z^{-1}] \ar[r]^{1-e+ze}
\ar[d]^g & P[z,z^{-1}] \ar[d]^g\ar[dl]_h \\
P'[z,z^{-1}] \ar[r]^{1-e'+ze'} & P'[z,z^{-1}]~.}$$
Thus
$$h(1-e+ze)~=~g~:~P[z,z^{-1}] \to P'[z,z^{-1}]~,$$
and writing
$$h~=~\sum\limits^b_{j=-a}z^jh_j~:~P[z,z^{-1}] \to P'[z,z^{-1}]$$
we have
$$h_{j-1}e+h_j(1-e)~=~
\begin{cases} g&\hbox{if $j=0$} \cr 0 & \hbox{if $j \neq 0$} \end{cases}~.$$
For any $k \geqslant 1$
$$\begin{array}{ll}
g(e(1-e))^k&=~h_{-1}e^{k+1}(1-e)^k+h_0e^k(1-e)^{k+1}\\[1ex]
&=~-h_{-2}e^{k+2}(1-e)^{k-1}-h_1e^{k-1}(1-e)^{k+2}\\[1ex]
&=~h_{-3}e^{k+3}(1-e)^{k-2}+h_2e^{k-2}(1-e)^{k+3}\\[1ex]
&=~\dots\\[1ex]
&=~(-1)^k(h_{-k-1}e^{2k+1}+h_k(1-e)^{2k+1})~.
\end{array}$$
Now $h_{-k-1}=0$ for $k \geqslant a$, and $h_k=0$ for $k \geqslant b+1$, so that
for $k={\rm max}(a,b+1)$ we have
$$g(e(1-e))^k~=~0~:~P \to P'~.$$
\end{proof}

\begin{exa} \label{Sexa}
{\rm  Let $p:M \to S^1$ be a map from a finite $CW$ complex which
is transverse regular at a point $1 \in S^1$ in the sense that
$N=p^{-1}(1) \subset M$ is a subcomplex, and cutting $M$ along $N$
gives a fundamental domain $(M_N;N,zN)$ for the pullback infinite
cyclic cover of $M$
$$\overline{M}~=~p^*\RR~=~\bigcup\limits^{\infty}_{j=-\infty}z^jM_N$$
with $z:\overline{M} \to \overline{M}$ a generating covering
translation. The map $p$ can be cut also, to obtain a map
$$p_N~:~(M_N;N,zN) \to ([0,1];\{0\},\{1\})$$
such that
$$p~=~[p_N]~:~M~=~M_N/(N=zN) \to S^1~=~[0,1]/(0=1)~.
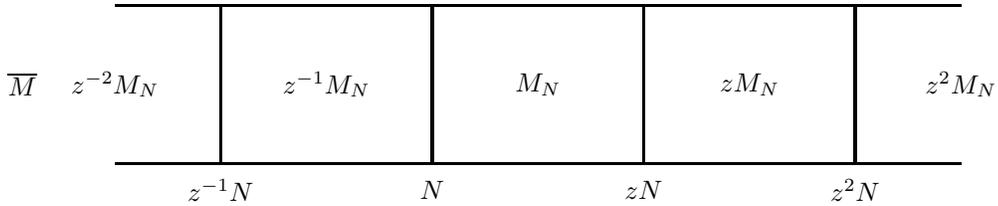
$$
$$\Draw
\LineAt(-160,-30,160,-30) \LineAt(-160,30,160,30)
\LineAt(-120,30,-120,-30) \LineAt(-40,30,-40,-30)
\LineAt(40,30,40,-30) \LineAt(120,30,120,-30) \MoveTo(-195,0)
\Text(--$\overline{M}$--) \MoveTo(-160,0) \Text(--$z^{-2}M_N$--)
\MoveTo(-80,0) \Text(--$z^{-1}M_N$--) \MoveTo(0,0)
\Text(--$M_N$--) \MoveTo(80,0) \Text(--$zM_N$--) \MoveTo(160,0)
\Text(--$z^2M_N$--) \MoveTo(-120,-40) \Text(--$z^{-1}N$--)
\MoveTo(-40,-40) \Text(--$N$--) \MoveTo(40,-40) \Text(--$zN$--)
    \MoveTo(120,-40) \Text(--$z^2N$--) \EndDraw$$
The two inclusions
$$f~:~N \to M_N~~,~~g~:~N~=~zN \to M_N$$
induce chain maps of finite f.g. free $\Z$-module chain complexes
$$f,g~:~C~=~C(M \to \{0\})_{*+1} \to D~=~C(p_N:M_N \to [0,1])_{*+1}$$
such that
    $$C(f-zg:C[z,z^{-1}] \to
    D[z,z^{-1}])~=~C(\overline{p}:\overline{M} \to \RR)_{*+1}~.$$
In particular, if $M$ is a knot complement then $p:M \to S^1$ can
be chosen to be a homology equivalence, and $N \subset M$ is a
Seifert surface for the knot, as in the Introduction and Example
\ref{Bexa}. In this case
$$H_*(f-g)~=~H_{*+1}(\overline{p}:\overline{M} \to \RR)~=~0$$
and $f-g:C \to D$ is a chain equivalence. The $\Z$-module chain
map
$$e~=~(f-g)^{-1}f~:~C \to C$$
defines a finite chain complex $(C,e)$ in the Seifert module
category $\E(\Z)$ with covering $B(C,e)$ a finite chain complex in
the Blanchfield module category $\B(\Z[z,z^{-1}])$ such that
$$B(C,e)~\simeq~C(\overline{p}:\overline{M} \to \RR)_{*+1}~.\eqno{\qed}$$
}
\end{exa}

\section{Seifert modules with zero Blanchfield module}\label{kernel}

This Chapter is devoted to the kernel of the covering functor
from Seifert modules to Blanchfield modules
$$B~:~\E(A) \to \B(A[z,z^{-1}])~;~(P,e) \mapsto B(P,e)~.$$
We study the Seifert modules $(P,e)$
with $B(P,e)=0$, and more generally the morphisms of Seifert
modules $g:(P,e) \to (P',e')$ with $B(g):B(P,e) \to B(P',e')$ an
isomorphism.

\begin{defn} {\rm
(i) A  Seifert $A$-module $(P,e)$ is {\it nilpotent} if
$$e^k~=~0~:~P \to P$$
for some $k \geqslant 0$.\\
(ii) A  Seifert $A$-module $(P,e)$ is {\it unipotent} if
$(P,1-e)$
is nilpotent, that is
$$(1-e)^k~=~0~:~P \to P$$
for some $k \geqslant 0$.\\
(iii) A  Seifert $A$-module $(P,e)$ is a {\it projection} if
$$e(1-e)~=~0~:~P \to P~.$$
(iv) A  Seifert $A$-module $(P,e)$ is a {\it near-projection} if
$e(1-e):P \to P$ is nilpotent, that is if for some $k \geqslant 0$
$$(e(1-e))^k~=~0~:~P \to P~.\eqno{\qed}$$}
\end{defn}

The near-projection terminology was introduced in L\"uck and Ranicki \cite{LR}.

\begin{prop} \label{bhs}
{\rm (Bass, Heller and Swan \cite{BHS})}\\
{\rm (i)} A linear morphism of induced f.g. projective $A[z]$-modules
$$f_0+zf_1~:~P[z] \to Q[z]$$
is an isomorphism if and only if $f_0+f_1:P \to Q$ is an isomorphism and
$$e~=~(f_0+f_1)^{-1}f_1~:~P \to P$$
is nilpotent.\\
{\rm (ii)} A linear morphism of induced f.g. projective $A[z,z^{-1}]$-modules
$$f_0+zf_1~:~P[z,z^{-1}] \to Q[z,z^{-1}]$$
is an isomorphism if and only if $f_0+f_1:P \to Q$ is an isomorphism and
$$e~=~(f_0+f_1)^{-1}f_1~:~P \to P$$
is a near-projection.\hfill\qed
\end{prop}

\begin{prop} \label{decomp0}
The following conditions on a Seifert $A$-module $(P,e)$ are equivalent :
\begin{itemize}
\item[(i)] $B(P,e)=0$.
\item[(ii)] $(P,e)$ is a near-projection.
\item[(iii)] There is a direct sum decomposition
$$(P,e)~=~(P^+,e^+) \oplus (P^-,e^-)$$
with $(P^+,e^+)$ unipotent and $(P^-,e^-)$ nilpotent.
\end{itemize}
\end{prop}
\begin{proof}
(i) $\Longleftrightarrow$ (ii) This is a special case of Proposition
\ref{bhs}, with
$$f~=~1-e+ze~:~P[z,z^{-1}] \to P[z,z^{-1}]~.$$
(iii) $\Longrightarrow$ (ii) Immediate from
$$e(1-e)~=~e^+(1-e^+) \oplus e^-(1-e^-)~:~P~=~P^+ \oplus P^- \to
P~=~P^+ \oplus P^-~.$$
(ii) $\Longrightarrow$ (iii) By the binomial theorem, for any $k \geqslant 1$ and an
indeterminate $x$ over ${\mathbb Z}$
$$x^k+(1-x)^k~=~1+x(1-x)\pi_k(x) \in {\mathbb Z}[x]$$
with
$$\pi_k(x)~=~\sum\limits^{k-1}_{j=1}((-)^j\binom{k-1}{j}-1)x^{j-1}
\in {\mathbb Z}[x]~.$$
Thus for any $A$-module endomorphism $e:P \to P$
$$e^k+(1-e)^k~=~1+e(1-e)\pi_k(e)~:~P \to P~.$$
If $(P,e)$ is a near-projection with $(e(1-e))^k=0$ then
$e(1-e)\pi_k(e):P \to P$ is nilpotent, and $e^k+(1-e)^k:P \to P$
is an automorphism. The endomorphism
$$p~=~(e^k+(1-e)^k)^{-1}e^k~:~P \to P$$
is a projection, $p^2=p$, and the images
$$P^+~=~{\rm im}(p:P \to P)~~,~~P^-~=~{\rm im}(1-p:P \to P)$$
are such that
$$(P,e)~=~(P^+,e^+) \oplus (P^-,e^-)$$
with
$$(1-e^+)^k~=~0~:~P^+ \to P^+~~,~~(e^-)^k~=~0~:~P^- \to P^-~.$$
\end{proof}

\begin{prop} \label{decomp1}
Given a morphism $g:(P_1,e_1)\to (P_0,e_0)$
of Seifert $A$-modules let  $C$ be the 1-dimensional f.g. projective
$A$-module chain complex
$$d_C~=~g~:~C_1~=~P_1 \to C_0~=~P_0$$
and let $e:C \to C$ be the $A$-module chain map defined by
$$\begin{array}{l}
e_0~:~C_0~=~P_0 \to C_0~=~P_0~,\\[1ex]
e_1~:~C_1~=~P_1 \to C_1~=~P_1~.
\end{array}$$
The following conditions on $g$ are equivalent :
\begin{itemize}
\item[(i)] $B(g):B(P_1,e_1) \to B(P_0,e_0)$ is an isomorphism
of Blanchfield $A[z,z^{-1}]$-modules.
\item[(ii)] There exists a morphism $h:(P_0,e_0) \to (P_1,e_1)$
of Seifert $A$-modules such that
$$\begin{array}{l}
gh~=~(e_0(1-e_0))^k~:~P_0 \to P_0~,\\[1ex]
hg~=~(e_1(1-e_1))^k~:~P_1 \to P_1
\end{array}$$
for some $k \geqslant 0$, defining a chain homotopy
$$h~:~(e(1-e))^k~ \simeq~ 0~:~C \to C~.$$
\item[(iii)] There exist 1-dimensional f.g. projective $A$-module chain
complexes $C^+,C^-$ with chain maps
$$e^+~:~C^+ \to C^+~~,~~e^-~:~C^- \to C^-$$
such that $1-e^+:C^+ \to C^+$, $e^-:C^- \to C^-$ are chain homotopy nilpotent,
and with a chain equivalence
$$i~=~\begin{pmatrix} i^+ \\ i^- \end{pmatrix}~:~C \to C^+ \oplus C^-$$
such that
$$e^+i^+~=~i^+e~:~C \to C^+~~,~~e^-i^-~=~i^-e~:~C \to C^-~.$$
\end{itemize}
\end{prop}
\begin{proof}
(i) $\Longrightarrow$ (ii) By Theorem \ref{cover1} (iii) there
exist a morphism $i:(P_0,e_0) \to (P_1,e_1)$ and $j\geqslant 0$
such that
$$B(g)^{-1}~=~B(i)t^{-j}~:~B(P_0,e_0) \to B(P_1,e_1)~.$$
It follows that
$$\begin{array}{l}
B(gi)~=~B(g)B(i)~=~t^{-j}~=~B((e_0(1-e_0))^j)~:~B(P_0,e_0) \to B(P_0,e_0)~,\\[1ex]
B(ig)~=~B(i)B(g)~=~t^{-j}~=~B((e_1(1-e_1))^j)~:~B(P_1,e_1) \to B(P_1,e_1)
\end{array}$$
and by Theorem \ref{cover1} (iv) there exist $\ell_0,\ell_1
\geqslant 0$ such that
$$\begin{array}{l}
(gi-(e_0(1-e_0))^j)(e_0(1-e_0))^{\ell_0}~=~0~:~(P,e_0) \to (P,e_0)~,\\[1ex]
(ig-(e_1(1-e_1))^j)(e_1(1-e_1))^{\ell_1}~=~0~:~(P_1,e_1) \to (P_1,e_1)~.
\end{array}$$
Defining
$$\begin{array}{l}
h~=~i(e_0(1-e_0))^{\ell_0+\ell_1}~:~(P_0,e_0) \to (P_1,e_1)~,\\[1ex]
k~=~j+\ell_0+\ell_1
\end{array}$$
we have
$$\begin{array}{l}
gh~=~gi(e_0(1-e_0))^{\ell_0+\ell_1}~=~(e_0(1-e_0))^k~:~(P_0,e_0) \to (P_0,e_0)~,\\[1ex]
hg~=~ig(e_1(1-e_1))^{\ell_0+\ell_1}~=~(e_1(1-e_1))^k~:~(P_1,e_1) \to (P_1,e_1)~.
\end{array}$$
(ii) $\Longrightarrow$ (i) The inverse of $B(g)$ is given by
$$B(g)^{-1}~=~B(h)t^{-k}~:~B(P_0,e_0) \to B(P_1,e_1)~.$$
(iii) $\Longrightarrow$ (i) It follows from the chain homotopy
nilpotence of $1-e^+$ and $e^-$ that the $A[z,z^{-1}]$-module chain maps
$$\begin{array}{l}
1-e^++ze^+~:~C^+[z,z^{-1}] \to C^+[z,z^{-1}]~,\\[1ex]
1-e^-+ze^-~:~C^-[z,z^{-1}] \to C^-[z,z^{-1}]
\end{array}$$
are chain equivalences. It now follows from the commutative diagram
$$\xymatrix{C[z,z^{-1}] \ar[r]^-i \ar[d]_-{1-e+ze} &
(C^+ \oplus C^-)[z,z^{-1}] \ar[d]^-{(1-e^++ze^+)\oplus (1-e^-+ze^-)} \\
C[z,z^{-1}] \ar[r]^-i & (C^+ \oplus C^-)[z,z^{-1}]}$$
that the $A[z,z^{-1}]$-module chain map
$$1-e+ze~:~C[z,z^{-1}] \to C[z,z^{-1}]$$
is also a chain equivalence. Thus
$$\begin{array}{l}
{\rm coker}(B(g):B(P_1,e_1) \to B(P_0,e_0))~=~H_0(1-e+ze)~=~0~,\\[1ex]
{\rm ker}(B(g):B(P_1,e_1) \to B(P_0,e_0))~=~H_1(1-e+ze)~=~0
\end{array}$$
and $B(g):B(P_1,e_1) \to B(P_0,e_0)$ is an isomorphism.\\
(ii) $\Longrightarrow$ (iii)
As in the proof of Proposition \ref{decomp0} write
$$x^k+(1-x)^k~=~1+x(1-x)\pi_k(x) \in {\mathbb Z}[x]~.$$
The $A$-module chain map
$$e^k+(1-e)^k~:~C \to C$$
is a chain equivalence, with the $A$-module morphisms
$$\begin{array}{l}
u_0~=~\sum\limits^{k-1}_{j=0}(-e_0(1-e_0)\pi_k(e_0))^j~:~P_0 \to P_0~,\\[1ex]
u_1~=~\sum\limits^{k-1}_{j=0}(-e_1(1-e_1)\pi_k(e_1))^j~:~P_1 \to P_1
\end{array}$$
defining a chain homotopy inverse $u:C \to C$, and
the $A$-module morphism
$$v~=~(-\pi_k(e_1))^kh~:~P_0 \to P_1$$
defining a chain homotopy
$$v~:~u(e^k+(1-e)^k)~\simeq~1~:~C \to C~.$$
The $A$-module chain map
$$p~=~ue^k~:~C \to C$$
is a chain homotopy projection, with a chain homotopy
$$v~:~u(1-e)^k~\simeq~1-p~:~C \to C~,$$
and the $A$-module morphism
$$q~=~(u_1)^2h+pv~=~\sum\limits^k_{j=0}(-e_1(1-e_1)\pi_k(e_1))^ju_1h~:~P_0 \to P_1$$
defining a chain homotopy
$$q~:~p(1-p)~\simeq~0~:~C \to C~.$$
The $A$-module morphisms
$$\begin{array}{l}
p^+~=~\begin{pmatrix} p_0 & g \\ q & 1-p_1 \end{pmatrix}~:~
P_0 \oplus P_1 \to P_0 \oplus P_1~,\\[2ex]
p^-~=~\begin{pmatrix} 1-p_0 & -g \\ -q & p_1 \end{pmatrix}~:~
P_0 \oplus P_1 \to P_0 \oplus P_1
\end{array}$$
are projections such that
$$p^+ + p^-~=~1~:~P_0 \oplus P_1 \to P_0 \oplus P_1~.$$
(This is a special case of the {\it instant finiteness obstruction}
of Ranicki \cite{RFO} and L\"uck and Ranicki \cite{LR}).
Define 1-dimensional f.g. projective $A$-module chain complexes
$C^+,C^-$ by
$$\begin{array}{l}
d_{C^+}~=~p^+\vert~:~C^+_1~=~P_1 \to C^+_0~=~{\rm im}(p^+)~,\\[1ex]
d_{C^-}~=~p^-\vert~:~C^-_1~=~P_1 \to C^-_0~=~{\rm im}(p^-)~.
\end{array}$$
The $A$-module chain maps
$$e^+~:~C^+ \to C^+~,~e^-~:~C^- \to C^-~,~i^+~:~C \to C^+~,~
i^-~:~C \to C^-$$
defined by
$$\begin{array}{l}
e^+_0~=~(e_0 \oplus e_1)\vert~:~C^+_0~=~{\rm im}(p^+) \to C^+_0~=~{\rm im}(p^+)~,\\[1ex]
e^+_1~=~e_1~:~C^+_1~=~P_1 \to C^+_1~=~P_1~,\\[1ex]
e^-_0~=~(e_0 \oplus e_1)\vert~:~C^-_0~=~{\rm im}(p^-) \to C^-_0~=~{\rm im}(p^-)~,\\[1ex]
e^-_1~=~e_1~:~C^-_1~=~P_1 \to C^-_1~=~P_1~,\\[1ex]
i^+_0~=~p^+\vert~:~C_0~=~P_0 \to C^+_0~=~{\rm im}(p^+)~,\\[1ex]
i^+_1~=~p_1~:~C_1~=~P_1 \to C^+_1~=~P_1~,\\[1ex]
i^-_0~=~p^-\vert~:~C_0~=~P_0 \to C^-_0~=~{\rm im}(p^-)~,\\[1ex]
i^-_1~=~1-p_1~:~C_1~=~P_1 \to C^-_1~=~P_1
\end{array}$$
are such that $1-e^+:C^+ \to C^+$, $e^-:C^- \to C^-$ are chain
homotopy nilpotent, with
$$i~=~\begin{pmatrix} i^+ \\ i^- \end{pmatrix}~:~C \to C^+ \oplus C^-$$
a chain equivalence such that
$$e^+i^+~=~i^+e~:~C \to C^+~~,~~e^-i^-~=~i^-e~:~C \to C^-~.$$
\end{proof}

\begin{rem} \label{decomp3} {\rm
(a) Propositions \ref{decomp0}
and \ref{decomp1} are the 0- and 1-dimensional cases of a general
result, namely that the following conditions on a self chain map
$e:C \to C$ of an $n$-dimensional
f.g. projective $A$-module chain complex are equivalent:\begin{itemize}
\item[(i)] The $A[z,z^{-1}]$-module chain map
$$1-e+ze~:~C[z,z^{-1}] \to C[z,z^{-1}]$$
is a chain equivalence.
\item[(ii)] For some $k \geqslant 0$ there exists a chain homotopy
$$h~:~(e(1-e))^k~\simeq~0~:~C \to C$$
such that $eh=he$, i.e. $e:C \to C$ is a chain homotopy near-projection.
\item[(iii)] There exist $n$-dimensional f.g. projective $A$-module chain
complexes $C^+,C^-$ with chain maps
$$e^+~:~C^+ \to C^+~~,~~e^-~:~C^- \to C^-$$
such that $1-e^+:C^+ \to C^+$, $e^-:C^- \to C^-$ are chain homotopy nilpotent,
and with a chain equivalence
$$i~=~\begin{pmatrix} i^+ \\ i^- \end{pmatrix}~:~C \to C^+ \oplus C^-$$
such that
$$e^+i^+~=~i^+e~:~C \to C^+~~,~~e^-i^-~=~i^-e~:~C \to C^-~.$$
\end{itemize}
(b) If $(C,e)$ satisfies the equivalent conditions in (a) then there there
are defined $A$-module chain equivalences
$$\begin{array}{l}
C(1-e+ze:C[z] \to C[z])~\simeq~C(1-e^++ze^+:C^+[z] \to C^+[z])~
\simeq~C^+~,\\[1ex]
C(z^{-1}(1-e)+e:C[z^{-1}] \to C[z^{-1}])~\simeq\\[1ex]
\hskip120pt C(z^{-1}(1-e^-)+e^-:C^-[z^{-1}] \to C^-[z^{-1}])~\simeq~C^-~,
\end{array}$$
so that the chain homotopy types of $C^+,C^-$ are entirely determined
by $C$ and $e$.
\hfill\qed}
\end{rem}

\section{Blanchfield and Seifert forms}\label{forms}

Let now $A$ be a ring with involution $A \to A;a \mapsto \overline{a}$.

\begin{defn} \label{dual}
{\rm
(i) The {\it dual} of a f.g. projective (left) $A$-module $P$ is
the f.g. projective $A$-module
$$P^*~=~{\rm Hom}_A(P,A)$$
with
$$A \times P^* \to P^*~;~(a,f) \mapsto (x \mapsto f(x)\overline{a})~.$$
(ii) The {\it dual} of a morphism $f:P \to Q$ of f.g. projective
$A$-modules is the morphism
$$f^*~:~Q^* \to P^*~;~g \mapsto (x \mapsto g(f(x)))~.\eqno{\qed}$$}
\end{defn}

The natural $A$-module morphism
$$P \to P^{**}~;~x \mapsto (f \mapsto \overline{f(x)})$$
is an isomorphism, which will be used to identify
$$P^{**}~=~P~.$$
Thus for any f.g. projective $A$-modules duality defines an isomorphism
$$T~:~{\rm Hom}_A(P,Q) \to {\rm Hom}_A(Q^*,P^*)~;~f \mapsto f^*$$
with inverse $g \mapsto g^*$. In particular, for $Q=P^*$ this is an involution
$$T~:~{\rm Hom}_A(P,P^*) \to {\rm Hom}_A(P,P^*)~;~f \mapsto f^*$$
with $T^2=1$.

Fix a central unit $\eta \in A$ such that
$$\overline{\eta}~=~\eta^{-1} \in A~.$$
In practice, $\eta=+1$ or $-1$.

\begin{defn} \label{symmetric form}
{\rm An {\it $\eta$-symmetric form over $A$} $(P,\lambda)$ is a f.g.
projective $A$-module $P$ together with a morphism $\lambda:P \to P^*$
such that
$$\eta \lambda^*~=~\lambda~:~P \to P^*~.$$
The form is {\it nonsingular} if $\lambda:P\to P^*$ is an isomorphism.\hfill\qed}
\end{defn}

Extend the involution on $A$ to an involution on $A[z,z^{-1}]$ by
$$\overline{z}~=~z^{-1}~.$$

\begin{defn} \label{Bdual1}
{\rm
(i) The {\it dual} of a Blanchfield $A[z,z^{-1}]$-module
$B$ is the Blanchfield $A[z,z^{-1}]$-module
$$B\widehat{~}~=~{\rm Ext}^1_{A[z,z^{-1}]}(B,A[z,z^{-1}])~.$$
(ii) The {\it dual} of a Seifert $A$-module $(P,e)$ is the
Seifert $A$-module
$$(P,e)^*~=~(P^*,1-e^*)~.\eqno{\qed}$$}
\end{defn}

\begin{prop} {\rm (i)}
The dual of an induced f.g.  projective $A[z,z^{-1}]$-module
presentation of a Blanchfield $A[z,z^{-1}]$-module $B$
$$C~:~0 \to P_1[z,z^{-1}] \xymatrix{\ar[r]^{d}&} P_0[z,z^{-1}] \to B \to 0$$
is an induced f.g.  projective $A[z,z^{-1}]$-module presentation of
the dual Blanchfield $A[z,z^{-1}]$-module $B\widehat{~}$
$$C^{1-*}~:~0 \to P^0[z,z^{-1}] \xymatrix{\ar[r]^{d^*}&} P^1[z,z^{-1}] \to B\widehat{~} \to 0$$
with $P^i=(P_i)^*$ the dual f.g. projective $A$-modules.\\
{\rm (ii)} The dual $B(P,e)\widehat{~}$ of the covering $B(P,e)$
of a Seifert $A$-module $(P,e)$ is related to the covering
$B((P,e)^*)=B(P^*,1-e^*)$ of the dual Seifert $A$-module by a
natural isomorphism
$$\zeta_{(P,e)}~:~B(P^*,1-e^*) \to B(P,e)\widehat{~}~.$$
{\rm (iii)} For any Blanchfield $A[z,z^{-1}]$-module $B$ there is a natural
isomorphism
$$B~\cong~B\widehat{~}\widehat{~}~.$$
\end{prop}
\begin{proof} (i) Any exact sequence of projective $A[z,z^{-1}]$-modules
$$0 \to Q_1 \to Q_0 \to B \to 0$$
induces an exact sequence
$$\begin{array}{l}
{\rm Hom}_{A[z,z^{-1}]}(B,A[z,z^{-1}])=0 \to
{\rm Hom}_{A[z,z^{-1}]}(Q_0,A[z,z^{-1}])\\[1ex]
\hskip50pt  \to {\rm Hom}_{A[z,z^{-1}]}(Q_1,A[z,z^{-1}])\\[1ex]
\hskip50pt \to{\rm Ext}^1_{A[z,z^{-1}]}(B,A[z,z^{-1}]) \to
{\rm Ext}^1_{A[z,z^{-1}]}(Q_0,A[z,z^{-1}])=0~.
\end{array}$$
(ii) Define $\zeta_{(P,e)}$ to fit into the natural isomorphism of exact sequences
of induced f.g. projective $A[z,z^{-1}]$-modules
$$\xymatrix@C+20pt{0\ar[r] &
P^*[z,z^{-1}] \ar[r]^{e^*+z(1-e^*)}
\ar[d]^z & P^*[z,z^{-1}] \ar[d]^1 \ar[r] & B(P^*,1-e^*)\ar[d]^{\zeta_{(P,e)}}
\ar[r] & 0\\
0 \ar[r] &P^*[z,z^{-1}] \ar[r]^{1-e^*+z^{-1}e^*} &
P^*[z,z^{-1}] \ar[r] & B(P,e)\widehat{~}
\ar[r]& 0}$$
(iii) The double dual $C^{**}$ of any induced f.g. projective $A[z,z^{-1}]$-module
resolution
$$C~:~0 \to C_1 \to C_0 \to B \to 0$$
is naturally isomorphic to $C$.
\end{proof}

\begin{defn} \label{Bdual2}
{\rm (i) The {\it dual} of a morphism $g:(P,e) \to (P',e')$ of
Seifert $A$-modules is the morphism
$$g^*~:~({P'}^*,1-{e'}^*) \to (P^*,1-e^*)~.$$
(ii) The {\it dual} of a morphism $f:B \to B'$ of
Blanchfield $A[z,z^{-1}]$-modules is the morphism
$$f\widehat{~}~=~(f_1^*,f_0^*)~:~{B'}\widehat{~} \to B\widehat{~}$$
with $f_0,f_1$ the components of any chain map $C \to C'$
of induced f.g. projective $A[z,z^{-1}]$-module chain complexes
resolving $f$
$$\xymatrix{
0 \ar[r] &C_1 \ar[r]^d \ar[d]^{f_1} & C_0 \ar[d]^{f_0} \ar[r]
& B \ar[r] \ar[d]^f  & 0 \\
0 \ar[r] & C'_1 \ar[r]^{d'} & C'_0 \ar[r] & B' \ar[r] &0~,}$$
so that $(f_1^*,f_0^*)$ resolves $f\widehat{~}$
$$\xymatrix{
0 \ar[r] &{C'}^0 \ar[r]^{{d'}^*} \ar[d]^{f^*_0} & {C'}^1 \ar[d]^{f_1^*} \ar[r]
& {B'}\widehat{~} \ar[r] \ar[d]^{f\widehat{~}}  & 0 \\
0 \ar[r] & C^0 \ar[r]^{d^*} & C^1 \ar[r] & B\widehat{~} \ar[r] &0~.}$$
\hfill\qed}
\end{defn}

\begin{prop} \label{B}
{\rm (i)} For any Blanchfield $A[z,z^{-1}]$-modules $B,B'$ duality
defines an isomorphism
$$T~:~{\rm Hom}_{A[z,z^{-1}]}(B,B') \to
{\rm Hom}_{A[z,z^{-1}]}({B'}\widehat{~},B\widehat{~})~;~f \mapsto
f\widehat{~}$$
with inverse $g \mapsto g\widehat{~}$. For $B'=B\widehat{~}$ this is
an involution
$$T~:~{\rm Hom}_{A[z,z^{-1}]}(B,B\widehat{~}) \to
{\rm Hom}_{A[z,z^{-1}]}(B,B\widehat{~})~;~f \mapsto f\widehat{~}$$
with $T^2=1$.\\
{\rm (ii)} The dual of a morphism of Blanchfield $A[z,z^{-1}]$-modules
$$f~=~B(g)t^{-k}~:~B(P,e) \to B(P',e')$$
is the morphism
$$f\widehat{~}~:~B(P',e')\widehat{~} \to B(P,e)\widehat{~}$$
such that
$$
(\zeta_{(P,e)})^{-1}f\widehat{~}\zeta_{(P',e')}~=~B(g^*)t^{-k}~:~
B({P'}^*,1-{e'}^*) \to B(P^*,1-e^*)~.$$
{\rm (iii)} For any Seifert $A$-module $(P,e)$ the dual of the isomorphism
of Blanchfield $A[z,z^{-1}]$-modules
$\zeta_{(P,e)}:B(P^*,1-e^*) \to B(P,e)\widehat{~}$ is the isomorphism
$$(\zeta_{(P,e)})\widehat{~}~=~z^{-1}\zeta_{(P^*,1-e^*)}~:~
B(P,e) \to B(P^*,1-e^*)\widehat{~}~.$$
{\rm (iv)} For any Seifert $A$-module $(P,e)$ the duality involution
$$T~:~{\rm Hom}_{A[z,z^{-1}]}(B(P,e),B(P,e)\widehat{~})
\to {\rm Hom}_{A[z,z^{-1}]}(B(P,e),B(P,e)\widehat{~})~;~f \mapsto
f\widehat{~}$$
corresponds under the isomorphism induced by
$\zeta_{(P,e)}:B(P^*,1-e^*)\to B(P,e)\widehat{~}$
$$\begin{array}{l}
\zeta_{(P,e)}~:~{\rm Hom}_{A[z,z^{-1}]}(B(P,e),B(P^*,1-e^*))
\to {\rm Hom}_{A[z,z^{-1}]}(B(P,e),B(P,e)\widehat{~})~;\\[1ex]
\hskip200pt B(\theta)t^{-k} \mapsto \zeta_{(P,e)}B(\theta)t^{-k}
\end{array}
$$
to the $z^{-1}$-duality involution
$$\begin{array}{l}
T_{z^{-1}}~:~{\rm Hom}_{A[z,z^{-1}]}(B(P,e),B(P^*,1-e^*))\\[1ex]
\hskip150pt
\to {\rm Hom}_{A[z,z^{-1}]}(B(P,e),B(P^*,1-e^*))~;\\[1ex]
\hskip150pt B(\theta)t^{-k} \mapsto z^{-1}B(\theta^*)t^{-k}~.
\end{array}$$
\end{prop}
\begin{proof} (i) By construction.\\
(ii) Applying Definition \ref{Bdual2} to the resolution of $f$
$$\xymatrix@C+25pt{
0 \ar[r] &P[z,z^{-1}] \ar[r]^{1-e+ze} \ar[d]^{gt^{-k}} &
P[z,z^{-1}] \ar[d]^{gt^{-k}} \ar[r]
& B(P,e) \ar[r] \ar[d]^f  & 0 \\
0 \ar[r] & P'[z,z^{-1}] \ar[r]^{1-e'+ze'} & P'[z,z^{-1}] \ar[r] & B(P',e') \ar[r] &0}$$
the identity $(\zeta_{(P,e)})^{-1}f\widehat{~}\zeta_{(P',e')}=B(g^*)t^{-k}$
is given by the composition of resolutions
$$\xymatrix@C+15pt{
0 \ar[r] &{P'}^*[z,z^{-1}] \ar[r]^{{e'}^*+z(1-{e'}^*)} \ar[d]^z &
{P'}^*[z,z^{-1}] \ar[d]^1 \ar[r]
& B({P'}^*,1-{e'}^*) \ar[r] \ar[d]^{\zeta_{(P',e')}}  & 0 \\
0 \ar[r] &{P'}^*[z,z^{-1}] \ar[r]^{1-{e'}^*+z^{-1}{e'}^*} \ar[d]^{g^*t^{-k}} &
{P'}^*[z,z^{-1}] \ar[d]^{g^*t^{-k}} \ar[r]
& B(P',e')\widehat{~} \ar[r] \ar[d]^{f\widehat{~}}  & 0 \\
0 \ar[r] &P^*[z,z^{-1}] \ar[r]^{1-e^*+z^{-1}e^*} \ar[d]^{z^{-1}} &
P^*[z,z^{-1}] \ar[d]^1 \ar[r]
& B(P,e)\widehat{~} \ar[r] \ar[d]^{\zeta_{(P,e)}^{-1}}  & 0 \\
0 \ar[r] & P^*[z,z^{-1}] \ar[r]^{e^*+z(1-e^*)} & P^*[z,z^{-1}] \ar[r] &
B(P^*,1-e^*) \ar[r] &0}$$
(iii) Consider the composition of resolutions
$$\xymatrix@C+16pt{
0 \ar[r] &P[z,z^{-1}] \ar[r]^{1-e+ze} \ar[d]^1 &
P[z,z^{-1}] \ar[d]^{z^{-1}} \ar[r]
& B(P,e) \ar[r] \ar[d]^{(\zeta_{(P,e)})\widehat{~}}  & 0 \\
0 \ar[r] &P[z,z^{-1}] \ar[d]^{z^{-1}} \ar[r]^{e+z^{-1}(1-e)}  &
P[z,z^{-1}] \ar[r] \ar[d]^1
& B(P^*,1-e^*)\widehat{~} \ar[r] \ar[d]^{(\zeta_{(P^*,1-e^*)})^{-1}} & 0 \\
0 \ar[r] &P[z,z^{-1}] \ar[r]^{1-e+ze} &
P[z,z^{-1}]  \ar[r]& B(P,e) \ar[r] & 0 }$$
(iv) By (ii) and (iii), for any morphism $\theta:(P,e) \to (P^*,1-e^*)$
$$(\zeta_{(P,e)}B(\theta))\widehat{~}~=~z^{-1}(\zeta_{(P,e)}B(\theta^*))~:~
\theta,e) \to B(P,e)\widehat{~}~.$$
\end{proof}

\begin{defn} \label{Bcover} {\rm
(i) An {\it $\eta$-symmetric Blanchfield form over
$A[z,z^{-1}]$} $(B,\phi)$ is a Blanchfield $A[z,z^{-1}]$-module $B$
together with a morphism $\phi:B \to B\widehat{~}$ such that
$$\eta \phi\widehat{~}~=~\phi ~:~B \to B\widehat{~}~.$$
The form is {\it nonsingular} if $\phi:B \to B\widehat{~}$ is an
isomorphism.\\
A {\it morphism} of Blanchfield forms
$f:(B,\phi) \to (B',\phi')$ is a morphism of Blanchfield modules
$f:B \to B'$ such that
$$f\widehat{~}\phi'f~=~\phi~:~B \to B\widehat{~}~.$$
(ii) A {\it $(-\eta)$-symmetric Seifert form over $A$} $(P,e,\theta)$
is a morphism of Seifert $A$-modules
$$\theta~:~(P,e) \to (P^*,1-e^*)$$
such that
$$\theta~=~(\theta-\eta\theta^*)e~:~P \to P^*~.$$
(This is equivalent to a  morphism of Seifert $A$-modules
 $\lambda:(P,e) \to (P^*,1-e^*)$ such that $\eta\lambda^*=-\lambda$,
with $\theta=\lambda e$, $\theta-\eta \theta^*=\lambda$.)
The form $(P,e,\theta)$ is {\it nonsingular}
if $\theta-\eta \theta^*:P \to P^*$ is an isomorphism.\\
A {\it morphism} of Seifert forms
$g:(P,e,\theta)\to (P',e',\theta')$ is a morphism of Seifert modules
$g:(P,e) \to (P',e')$ such that
$$g^*\theta'g~=~\theta~:~P \to P^*~.$$
(iii) The {\it covering} of a $(-\eta)$-symmetric Seifert form
over $A$ $(P,e,\theta)$ is the $\eta$-symmetric Blanchfield form over $A[z,z^{-1}]$
$$B(P,e,\theta)~=~(B(P,e),\phi)$$
with
$$\phi~=~(1-z^{-1})\zeta_{(P,e)}B(\theta-\eta\theta^*)~:~B(P,e) \to B(P,e)\widehat{~}~.$$
If $(P,e,\theta)$ is a nonsingular Seifert form then $B(P,e,\theta)$ is
a nonsingular Blanchfield form.
\hfill\qed}
\end{defn}

\begin{exa} {\rm An $n$-knot $k:S^n \subset S^{n+2}$ with exterior $M$
determines a $\Z$-acyclic $(n+2)$-dimensional symmetric Poincar\'e
complex $(C,\phi)$ over $\Z[z,z^{-1}]$ with
$C=C(\overline{p}:\overline{M} \to \RR)_{*+1}$.  Furthermore, a Seifert
surface $N^{n+1} \subset S^{n+2}$ for $k$ determines an
$(n+1)$-dimensional Seifert $\Z$-module chain complex $(D,e,\theta)$
for $(C,\phi)$ with $D=\widetilde{C}(M)$ and $(C,\phi)=B(D,\theta)$.
If $n=2i-1$ and $\overline{M}$ is $(i-1)$-connected then $N$ can be
chosen to be $(i-1)$-connected; in this simple case $(H_i(C),\phi_0)$
is a nonsingular $(-1)^{i+1}$-symmetric Blanchfield form over
$\Z[z,z^{-1}]$, and $(H_i(D),e,\theta)$ is a nonsingular
$(-1)^i$-symmetric Seifert form over $\Z$ such that
$(H_i(C),\phi_0)=B(H_i(D),e,\theta)$, with
$e=(\theta+(-1)^i\theta^*)^{-1}\theta$.  See Ranicki \cite[Chapter
7.9]{RESATS}, \cite[Chapter 32]{RHK} for further details.
\hfill$\qed$}
\end{exa}

\begin{prop} \label{replace}
{\rm (i)} For any morphism from a Seifert $A$-module to its dual
$$\theta~:~(P,e) \to (P^*,1-e^*)$$
the morphism
$$\theta'~=~(\theta-\eta \theta^*)e~:~(P,e) \to (P^*,1-e^*)$$
defines a $(-\eta)$-symmetric Seifert form $(P,e,\theta')$ such that
$$\theta'-\eta{\theta'}^*~=~\theta-\eta\theta^*~:~P \to P^*~.$$
{\rm (ii)} For a nonsingular $(-\eta)$-symmetric Seifert form
$(P,e,\theta)$ the endomorphism $e:P \to P$ is determined by $\theta:P
\to P^*$, with
$$e~=~(\theta-\eta \theta^*)^{-1}\theta~:~P \to P~.$$
A morphism of nonsingular $(-\eta)$-symmetric Seifert forms
$g:(P,e,\theta) \to (P,e',\theta')$ is the same as a morphism
of the underlying Seifert modules $g:(P,e) \to (P',e')$ such that
$$g^*(\theta'-\eta\theta'^*)g~=~\theta-\eta\theta^*~:~P \to P^*~.$$
{\rm (iii)} Every morphism $f:B(P,e,\theta) \to B(P',e',\theta')$ of
the covering $\eta$-symmetric Blanchfield forms of $(-\eta)$-symmetric
Seifert forms $(P,e,\theta)$, $(P',e',\theta')$ is of the type
$$f~=~B(g)t^{-k}$$
with $k \geqslant 0$, $t=B(e(1-e)):B(P,e) \to B(P,e)$, and $g:(P,e) \to
(P',e')$ a morphism of Seifert $A$-modules such that for some $\ell
\geqslant 0$
    $$g^*(\theta'-\eta\theta'^*)g~=~
    (\theta-\eta\theta^*)(e(1-e))^{2\ell}~:~P \to P^*~.$$
\end{prop}
\begin{proof} (i) From the definitions
$$\begin{array}{ll}
\theta'-\eta{\theta'}^*
&=~(\theta-\eta\theta^*)e -\eta e^*(\theta^*-\overline{\eta}\theta)\\[1ex]
&=~(\theta-\eta\theta^*)e +(\theta-\eta\theta^*)(1-e)\\[1ex]
&=~\theta - \eta \theta^*~:~P \to P^*
\end{array}$$
and also
$$(\theta'-\eta{\theta'}^*)e~=~(\theta-\eta\theta^*)e~=~\theta'~:~P \to P^*~.$$
(ii) Immediate from the definitions.\\
(iii) By Theorem \ref{cover1} (iii) $f=B(h)t^{-j}$ for some $h:(P,e) \to (P',e')$, 
$j \geqslant 0$. Let $B(P,e,\theta)=(B(P,e),\phi)$, $B(P',e',\theta')=(B(P',e'),\phi')$,
so that $f\widehat{~}\phi'f=\phi$ and by  Proposition \ref{B} (ii) there
is defined a commutative diagram
    $$\xymatrix@C+15pt{B(P,e)  \ar[dr]^-{B(\theta-\eta\theta^*)}
    \ar[rrr]^{f=B(h)t^{-j}} \ar[dd]_-{\phi} &&&B(P',e')
    \ar[dl]_-{B(\theta'-\eta\theta')} \ar[dd]^-{\phi'} \\
    & B(P^*,1-e^*) \ar[dl]_-{(1-z^{-1})\zeta_{(P,e)}}
     & B({P'}^*,1-{e'}^*) \ar[dr]^{(1-z^{-1})\zeta_{(P',e')}}
    \ar[l]_-{B(h^*)t^{-j}} & \\
    B(P,e)\widehat{~} & & &\ar[lll]_-{f\widehat{~}}B(P',e')\widehat{~}}$$
Now apply \ref{cover1} (iv) to the identity
    $$B(h^*(\theta'-\eta\theta'^*)h)t^{-2j}~=~B(\theta-\eta\theta^*)~:~
      B(P,e) \to B(P^*,1-e^*)~,$$
to obtain 
$$(h^*(\theta'-\eta\theta'^*)h - (\theta-\eta\theta^*)(e(1-e))^{2j})
(e(1-e))^{\ell}~=~0~:~ B(P,e) \to B(P^*,1-e^*)$$
for some $\ell \geqslant 0$. Setting 
$$g~=~h(e(1-e))^\ell~,~k~=~j+\ell$$
gives the required expression $f=B(g)t^{-k}$.
\end{proof}

In Theorem \ref{cover2} below, it will be proved that every nonsingular
Blanchfield form over $A[z,z^{-1}]$ is isomorphic to the covering of
a nonsingular Seifert form over $A$. The proof will use the quadratic
Poincar\'e complexes of Ranicki \cite{RESATS}, \cite{RHK}.
By definition, a {\it 1-dimensional $\eta$-quadratic Poincar\'e complex}
$(C,\psi)$ over $A$ is a 1-dimensional f.g. projective $A$-module chain complex
$$\xymatrix@C-5pt{C~:~\dots \ar[r]& 0 \ar[r]&  C_1 \ar[r]^-{d} & C_0}$$
together with $A$-module morphisms
$$\psi_0~:~C^0~=~C^*_0 \to C_1~,~\widetilde{\psi}_0~:~C^1 \to C_0~,~
\psi_1~:~C^0 \to C_0$$
such that
$$d\psi_0+\widetilde{\psi}_0d^*+\psi_1-\eta\psi_1^*~=~0:~C^0 \to C_0$$
and the chain map
$(\psi_0+\eta\widetilde{\psi}^*_0,\widetilde{\psi}_0+\eta\psi_0^*):C^{1-*} \to C$
is a chain equivalence. Replacing $\psi_0,\widetilde{\psi}_0,\psi_1$ by
$\psi_0+\eta\widetilde{\psi}^*_0,0,\psi_1+\widetilde{\psi}_0d^*$
respectively, it may always be assumed that $\widetilde{\psi}_0=0$.

\begin{thm} \label{cover2}
Every nonsingular $\eta$-symmetric Blanchfield form
$(B,\phi)$ over $A[z,z^{-1}]$ is isomorphic to the covering $B(P,e,\theta)$
of a nonsingular $(-\eta)$-symmetric Seifert form $(P,e,\theta)$ over $A$.
If $B$ admits an induced f.g. projective $A[z,z^{-1}]$-module
resolution
$$0 \to P_1[z,z^{-1}] \xymatrix{\ar[r]^d&} P_0[z,z^{-1}] \to B \to 0$$
with
$$d~=~\sum\limits^k_{j=0}d_jz^j~:~P_1[z,z^{-1}] \to P_0[z,z^{-1}]$$
then $P$ can be chosen to be a direct summand of
$\bigoplus\limits_k(P_0\oplus P_0^*)$ such that
$$P\oplus P^*~\cong~\bigoplus\limits_k(P_0\oplus P_0^*)~.$$
\end{thm}
\begin{proof} By Proposition \ref{Bprop} the given resolution of $B$
determines a resolution of the form
$$0 \to P[z,z^{-1}] \xymatrix{\ar[r]^{1-e+ze}&} P[z,z^{-1}] \to B \to 0$$
with $e:P \to P$ an endomorphism of
$$P~=~\bigoplus\limits_kP_0~.$$
(This is not yet the $(P,e)$ we are seeking). By Theorem
\ref{cover1} (i), (iv) it may be assumed that
$$(\zeta_{(P,e)})^{-1}\phi t~=~B(\theta)t^{-\ell}~:~B~=~B(P,e) \to B(P^*,1-e^*)$$
for some Seifert $A$-module $(P,e)$, morphism
$\theta:(P,e) \to (P^*,1-e^*)$ and $\ell \geqslant 0$.
The $\eta$-symmetric Blanchfield form $(B(P,e),\phi')$ defined by
$$\phi'~=~\zeta_{(P,e)}B(\theta)t^{-1}~:~B(P,e) \to B(P,e)\widehat{~}$$
is nonsingular, and such that there is defined an isomorphism
$$s^{\ell}~:~(B(P,e),\phi') \to (B(P,e),\phi)~.$$
Replacing $(B(P,e),\phi)$ by $(B(P,e),\phi')$ it may thus be assumed that $\ell=0$,
with
$$(\zeta_{(P,e)})^{-1}\phi t~=~B(\theta)~:~B(P,e) \to B(P^*,1-e^*)~.$$
The covering of $\theta-\eta \theta^*:(P,e) \to (P^*,1-e^*)$ is the
isomorphism of Blanchfield $A[z,z^{-1}]$-modules
$$\begin{array}{ll}
B(\theta-\eta \theta^*)&=~(\zeta_{(P,e)})^{-1}\phi t -\eta z
(\zeta_{(P,e)})^{-1}\phi\widehat{~}t \\[1ex]
&=~(\zeta_{(P,e)})^{-1}(1-z)\phi t\\[1ex]
&=~(1-z^{-1})^{-1}(\zeta_{(P,e)})^{-1}\phi~:~B(P,e) \to B(P,e)\widehat{~}~.
\end{array}$$
Replacing $\theta$ by $\theta'=(\theta-\eta\theta^*)e$ (as in
Proposition \ref{replace} (i))
we have a $(-\eta)$-symmetric Seifert form $(P,e,\theta)$ such that
$$B(P,e,\theta)~\cong~(B,\phi)~.$$
However, in general $(P,e,\theta)$ may be singular, i.e.
$\theta-\eta\theta^*:P \to P^*$ need not be an isomorphism. We
shall obtain a nonsingular $(-\eta)$-symmetric Seifert form
$(P',e',\theta')$ such that $B(P',e',\theta') \cong (B,\phi)$ by
gluing together two null-cobordism of the 1-dimensional
$(-\eta)$-quadratic Poincar\'e complex $(C,\psi)$ defined by
$$\begin{array}{l}
d_C~=~\theta-\eta\theta^*~:~C_1~=~P \to C_0~=~P^*~,\\[1ex]
\psi_0~=~1~:~C^0~=~P \to C_1=~P~,\\[1ex]
\psi_1~=~-\theta~:~C^0~=~P \to C_0~=~P^*~.
\end{array}$$
One null-cobordism is easy: it is $(f:C \to D,(0,\psi))$ with
$$f~=~1~:~C_1~=~P \to D_1~=~P~~,~~D_i~=~0~{\rm for}~i \neq 1~.$$
The other null-cobordism is of the form
$(i^-:C \to C^-,(\delta\psi,\psi))$, with $i^-:C \to C^-$ constructed
by the method of Remark \ref{decomp3},
as follows. By Proposition \ref{decomp1} (ii) (with $g=\theta-\eta\theta^*$) there exists
a morphism
$$h~:~(P^*,1-e^*) \to (P,e)$$
such that
$$\begin{array}{l}
h(\theta-\eta\theta^*)~=~(e(1-e))^k~:~P \to P~,\\[1ex]
(\theta-\eta\theta^*)h~=~(e^*(1-e^*))^k~:~P^* \to P^*
\end{array}$$
for some $k \geqslant 0$. Let $E:C \to C$ be the chain map defined by
$$\begin{array}{l}
E_0~=~1-e^*~:~C_0~=~P^* \to C_0~=~P^*~,\\[1ex]
E_1~=~e~:~C_1~=~P \to C_1~=~P~.
\end{array}$$
As in the proof of \ref{decomp1} (ii) $\Longrightarrow$
(iii) $h$ determines a chain homotopy projection
$$\begin{array}{ll}
p&=~(E^k+(1-E)^k)^{-1}E^k\\[1ex]
&=~\big(\sum\limits^{k-1}_{j=0}(-E(1-E)\pi_k(E))^j\big)E^k:~C \to C
\end{array}$$
with a chain homotopy
$$q~:~p(1-p)~\simeq~0~:~C \to C$$
such that $pE=Ep$, $Eq=qE$, and such that
$$\begin{array}{l}
p^+~=~\begin{pmatrix} p_0 & \theta-\eta\theta^* \\ q & 1-p_1 \end{pmatrix}~:~
P^* \oplus P \to P^* \oplus P~,\\[2ex]
p^-~=~\begin{pmatrix} 1-p_0 & -(\theta-\eta\theta^*) \\ -q & p_1 \end{pmatrix}~:~
P^* \oplus P \to P^* \oplus P
\end{array}$$
are projections with
$$p^+ + p^-~=~1~:~P^* \oplus P \to P^* \oplus P~.$$
We now have a decomposition of Seifert $A$-modules
$$(P^*\oplus P,(1-e^*)\oplus e)~=~(P^+,e^+)\oplus (P^-,e^-)$$
with
$$P^+~=~{\rm im}(p^+)~,~P^-~=~{\rm im}(p^-)~.$$
The  1-dimensional f.g. projective $A$-module chain complexes
$C^+,C^-$ defined by
$$\begin{array}{l}
d_{C^+}~=~p^+\vert~:~C^+_1~=~P \to C^+_0~=~P^+~,\\[1ex]
d_{C^-}~=~p^-\vert~:~C^-_1~=~P \to C^-_0~=~P^-
\end{array}$$
and the $A$-module chain maps
$$E^+~:~C^+ \to C^+~,~E^-~:~C^- \to C^-~,~i^+~:~C \to C^+~,~
i^-~:~C \to C^-$$
defined by
$$\begin{array}{l}
E^+_0~=~e^+~:~C^+_0~=~P^+ \to C^+_0~=~P^+~,\\[1ex]
E^+_1~=~e~:~C^+_1~=~P \to C^+_1~=~P~,\\[1ex]
E^-_0~=~e^-~:~C^-_0~=~P^- \to C^-_0~=~P^-~,\\[1ex]
E^-_1~=~1-e~:~C^-_1~=~P \to C^-_1~=~P~,\\[1ex]
i^+_0~=~p^+\vert~:~C_0~=~P^* \to C^+_0~=~P^+~,\\[1ex]
i^+_1~=~p_1~:~C_1~=~P \to C^+_1~=~P~,\\[1ex]
i^-_0~=~p^-\vert~:~C_0~=~P \to C^-_0~=~P^-~,\\[1ex]
i^-_1~=~1-p_1~:~C_1~=~P \to C^-_1~=~P
\end{array}$$
are such that $1-E^+:C^+ \to C^+$, $E^-:C^- \to C^-$ are chain
homotopy nilpotent, with
$$i~=~\begin{pmatrix} i^+ \\ i^- \end{pmatrix}~:~C \to C^+ \oplus C^-$$
a chain equivalence such that
$$E^+i^+~=~i^+E~:~C \to C^+~~,~~E^-i^-~=~i^-E~:~C \to C^-~.$$
Moreover, it follows from
$$E^*~=~E~:~C^{1-*}~=~C \to C^{1-*}~=~C$$
that
$$\begin{array}{l}
p^*~=~p~:~C^{1-*}~=~C \to C^{1-*}~=~C~,\\[1ex]
p_1~=~1-p_0^*~:~P \to P~.
\end{array}$$
The morphism
$$h'~=~eh-\eta h^*e^*~:~(P^*,1-e^*) \to (P,e)$$
is such that
$$\begin{array}{l}
h'(\theta-\eta\theta^*)~=~(e(1-e))^k~:~P \to P~,\\[1ex]
(\theta-\eta\theta^*)h'~=~(e^*(1-e^*))^k~:~P^* \to P^*
\end{array}$$
and replacing $h$ by $h'$ in the construction of $q$ gives a chain homotopy
$$q'~:~p(1-p)~\simeq~0~:~C \to C$$
such that
$$q'~=~\theta -\eta \theta^*~:~P^* \to P$$
with $\theta$ (resp. $-\eta\theta^*$) the contribution of $eh$ (resp.
$-\eta h^*e^*$), and $e\theta=\theta(1-e^*)$. The morphism of
Seifert $A$-modules defined by
$$\lambda~=~\begin{pmatrix}
\theta -\eta\theta^* & -\eta p_0^* \\
p_0 & \theta-\eta \theta^* \end{pmatrix}~:~
(P^* \oplus P,\begin{pmatrix}
1-e^* & 0 \\
0 & e \end{pmatrix}) \to (P \oplus P^*,\begin{pmatrix}
e & 0 \\
0 & 1-e^* \end{pmatrix})$$
is such that $(-\eta)\lambda^*=\lambda$, and restricts to an isomorphism
$$\lambda^+~:~(P^+,e^+) \to ((P^+)^*,1-(e^+)^*)~,$$
identifying $(P^+)^*={\rm im}((p^+)^*)$.
The $(-\eta)$-symmetric Seifert form over $A$ defined by
$$(P',e',\theta')~=~(P^+,e^+,\lambda^+e^+)$$
is nonsingular and such that
$$B(P',e',\theta')~\cong~(B,\phi)~.$$
\vskip-5mm
\end{proof}

\begin{rem}{\rm (i) The proof of Theorem \ref{cover2} minimizes the use of the
theory of algebraic Poincar\'e complexes. However, it is based on
an idea of infinite gluing which really is best expressed in this
language, specifically the quadratic $Q$-groups of an $A$-module
chain complex $C$
$$Q_n(C)~=~H_n(\Z_2;C\otimes_AC)$$
which are the central objects of the theory, with the generator $T
\in \Z_2$ acting by
$$T~:~C_p \otimes_A C_q \to C_q \otimes_A C_p~;~x\otimes y \mapsto
(-1)^{pq}y \otimes x~.$$ (There is a brief review in Chapter 20 of
\cite{RHK}). A chain map $f:C \to D$ induces morphisms in the
$Q$-groups
$$f_{\%}~:~Q_n(C) \to Q_n(D)$$
which depend only on the chain homotopy class of $f$. As in
Definition 24.1 of \cite{RHK}, given chain maps $f,g:C \to D$ let
$Q_*(f,g)$ be the relative $Q$-groups which fit into the exact
sequence
$$\dots \to Q_{n+1}(f,g) \to Q_n(C) \xymatrix{\ar[r]^{f_{\%}-g_{\%}}&}
Q_n(D) \to Q_n(f,g) \to \dots$$ and define a union operation
$$U~:~Q_n(f,g) \to Q_n(U(f,g))$$
with
$$U(f,g)~=~C(f-zg:C[z,z^{-1}] \to D[z,z^{-1}])$$
an $A[z,z^{-1}]$-module chain complex. An element
$(\delta\theta,\theta) \in Q_{n+1}(f,g)$ is an $(n+1)$-dimensional
quadratic pair over $A$
$$x~=~((f~g):C\oplus C \to D,(\delta\theta,\theta \oplus -\theta))$$
and the union is an $(n+1)$-dimensional quadratic complex over
$A[z,z^{-1}]$
$$U(x)~=~(U(f,g),U(\delta\theta,\theta))~.$$
The construction mimics the construction of an infinite cyclic
cover by gluing together $\Z$ copies of a fundamental domain. If
$x$ is a Poincar\'e pair then $U(x)$ is a Poincar\'e complex. The
chain complex ingredient in the proof of Theorem \ref{cover2} is
the following characterization of the pairs $x$ such that the
union $U(x)$ is contractible, i.e. such that
$$f-zg~:~C[z,z^{-1}] \to D[z,z^{-1}]$$
is a chain equivalence. This is the case if and only if $f-g:C \to
D$ is a chain equivalence and $(f-g)^{-1}f:C \to C$ is a chain
homotopy near-projection. Thus there is no loss of generality in
taking
$$C~=~D~,~f~=~1-e~,~g~=~-e$$
for a chain homotopy near-projection $e:C \to C$, and as in Remark
\ref{decomp3} there is defined a chain equivalence
$$i~:~(C,e) \to (C^+,e^+)\oplus (C^-,e^-)$$
with $1-e^+:C^+ \to C^+$, $e^-:C^-\to C^-$ chain homotopy
nilpotent. The background to the proof of Theorem \ref{cover2} is
the computation
$$Q_*(f,g)~=~H_*((1-e^+)\otimes e^-:C^+ \otimes_A C^- \to C^+ \otimes_A C^-)$$
so that an element $(\delta\theta,\theta) \in Q_{n+1}(f,g)$ is
determined by a chain map
$$\theta~:~(C^+)^{n-*} \to C^-$$
together with a chain homotopy
$$\delta \theta ~:~e^-\theta~\simeq~\theta(1-e^+)^*~:~(C^+)^{n-*} \to C^-~.$$
The quadratic pair $((f~g):C \oplus C \to D,(\delta\theta,\theta\oplus
-\theta))$
is Poincar\'e if and only if $\theta$ is a chain equivalence.\\
(ii)  Here is a geometric interpretation of (i). Let $X$ be a
finite $n$-dimensional Poincar\'e complex, and let $F:M \to X
\times S^1$ be a homotopy equivalence from a closed
$(n+1)$-dimensional manifold $M$. The restriction of $F$ to a
transverse inverse image is an $n$-dimensional normal map
$$G~=~F\vert~:~ N~=~F^{-1}(X \times \{*\}) \to X$$
and cutting $M$ along $N$ gives a fundamental domain for $F^*(X
\times \RR)=\overline{M}$ with a normal map
$$G_N~=~\overline{F}\vert~:~(M_N;N,zN) \to X \times ([0,1];\{0\},\{1\})$$
such that
$$\overline{F}~=~\bigcup\limits^{\infty}_{j=-\infty}z^jG_N~:~
\overline{M}~=~\bigcup\limits^{\infty}_{j=-\infty}z^jM_N \to X
\times \RR$$ is a $\Z$-equivariant lift of $F$.\\
\indent Define the kernel $\Z[\pi_1(X)]$-module chain complexes
$$\begin{array}{l}
C~=~C(G:C(N) \to C(X))_{*+1}~,\\[1ex]
D~=~C(G_N:C(M_N) \to C(X\times [0,1]))_{*+1}~.
\end{array}$$
The chain maps $i_0,i_1:C \to D$ induced by the inclusions $N \to
M_N$, $zN \to M_N$ are such that $f-zg:C[z,z^{-1}] \to
D[z,z^{-1}]$ is a chain equivalence, since $F:M \to X \times S^1$
is a homotopy equivalence. The infinitely generated free
$\Z[\pi_1(X)]$-module chain complexes
$$\begin{array}{l}
C^+~=~C(\overline{G}^+:\overline{M}^+ \to X \times \RR^+)_{*+1}~,\\[1ex]
C^-~=~C(\overline{G}^-:\overline{M}^- \to X \times \RR^-)_{*+1}
\end{array}$$
with
$$\overline{M}^+~=~\bigcup\limits^{\infty}_{j=0}z^jM_N~~,~~
\overline{M}^-~=~ \bigcup\limits^{-1}_{j=-\infty}z^jM_N$$ are such
that there is defined an exact sequence
$$0 \to C \to C^+ \oplus C^- \to C(\overline{F})_{*+1} \to 0$$
with $C(\overline{F})$ the algebraic mapping cone of the
$\Z[\pi_1(X)][z,z^{-1}]$-module chain equivalence
$\overline{F}:C(\overline{M}) \to C(X \times \RR)$. Thus there is
defined a $\Z[\pi_1(X)]$-module chain equivalence
$$C~ \simeq~ C^+ \oplus C^-~,$$
and $C^+,C^-$ are chain equivalent to finite f.g. projective
$\Z[\pi_1(X)]$-module chain complexes. The quadratic Poincar\'e
kernel of $G_N$ is determined as in (i) by a chain equivalence
$\theta:(C^+)^{n-*} \to C^-$.\\
\indent In particular, if $n=2i$ and $G,G_N$ are $i$-connected then
$$K_i(N)~=~H_i(C)~=~H_i(C^+) \oplus H_i(C^-)~
=~H_{i+1}(\overline{M}^+,N)\oplus H_{i+1}(\overline{M}^-,N)$$
with an isomorphism $\theta:H_i(C^+)^*\to H_i(C^-)$.
Every homology class in $K_i(N)$ is a sum of a class which dies on
the right and one which dies on the left; the reduced projective class
$$[H_i(C^+)]~=~-[H_i(C^-)] \in\widetilde{K}_0(\Z[\pi_1(X)])~=~
K_0(\Z[\pi_1(X)])/K_0(\Z)$$
is the obstruction to finding a basis of classes which all die on the
left (or all die on the right).  The reduced nilpotent projective class
    $$[H_i(C^+),1-e^+]~=~-[H_i(C^-),e^-] \in
    \widetilde{\rm Nil}_0(\Z[\pi_1(X)])~=~ {\rm
    Nil}_0(\Z[\pi_1(X)])/K_0(\Z)$$
is the Farrell-Hsiang \cite{FH} splitting obstruction of $F$, which is
0 if (and for $i \geqslant 3$ only if) $G:N \to X$ can be chosen to be
a homotopy equivalence, or equivalently $(M_N;N,zN)$ can be chosen to
be an $h$-cobordism.  The surgery obstruction
$$\begin{array}{ll}
\sigma(G)~=~(K_i(N),\lambda,\mu) &
\in {\rm ker}(L^h_{2i}(\Z[\pi_1(X)]) \to L^p_{2i}(\Z[\pi_1(X)]))\\[1ex]
&\hskip10pt
=~{\rm im}(\widehat{H}^{2i+1}(\Z_2;\widetilde{K}_0(\Z[\pi_1(X)]))
\to L^h_{2i}(\Z[\pi_1(X)]))
\end{array}$$
is represented by the hyperbolic $(-1)^i$-quadratic form on the
f.g. projective $\Z[\pi_1(X)]$-module $H_i(C^+)$, with
$$\begin{array}{l}
\lambda~=~\begin{pmatrix} 0 & \theta^{-1} \\
(-1)^i(\theta^*)^{-1} & 0\end{pmatrix}~:\\[1ex]
K_i(N)~=~H_i(C^+) \oplus H_i(C^-) \to
K_i(N)^*~=~H_i(C^+)^* \oplus H_i(C^-)^*~,\\[1ex]
\mu~:~K_i(N) \to \Z[\pi_1(X)]/\{x-(-1)^i\overline{x}\}~;~
(a^+,a^-) \mapsto \theta^{-1}(a^-)(a^+)~.
\end{array}$$
However, in Theorem \ref{cover2} it is the other case $n=2i+1$
which occurs.\hfill\qed}
\end{rem}

\section{Witt groups}\label{Witt}

This Chapter extends the results of Chapter \ref{forms} to the
algebraic $L$-groups of Blanchfield and Seifert forms, using
the algebraic theory of surgery (Ranicki \cite{RESATS}, \cite{RHK}).

Cohn \cite{C} constructed the universal localization $\sigma^{-1}R$ of a
ring $R$ inverting a set $\sigma$ of square matrices over $R$.
The canonical ring morphism $R \to \sigma^{-1}R$
is universally $\sigma$-inverting~: for any ring morphism
$f:R \to S$ such that $f(s)$ is invertible for every $s \in \sigma$
there is a unique ring morphism $\sigma^{-1}R \to S$ such that
$$f~:~R \to \sigma^{-1}R \to S~.$$
\indent
See \cite{RESATS} or \cite{RHK} for the expression of the free Wall quadratic
$L$-groups $L^h_n(R)=L_n(R)$ of a ring with involution $R$ as the
cobordism groups of $n$-dimensional quadratic Poincar\'e complexes
$(C,\psi)$ over $R$ with $C$ f.g.  free.  In particular, $L_{2i}(R)$ is
the Witt group of nonsingular $(-1)^i$-quadratic forms over $R$.

For an injective universal localization $R \to \sigma^{-1}R$
of rings with involution the quadratic $L$-groups
of $R$ and $\sigma^{-1}R$ are related by the exact sequence of
of Vogel \cite{V} and Neeman and Ranicki \cite{NR}
    $$\xymatrix{\dots \ar[r] &L_n(R) \ar[r] & L_n(\sigma^{-1}R)
    \ar[r]^{\partial} & L_n(R,\sigma) \ar[r] & L_{n-1}(R) \ar[r] & \dots}$$
with $L_n(R,\sigma)$ the cobordism group of $(n-1)$-dimensional
quadratic Poincar\'e complexes $(C,\psi)$ over $R$ such that $C$ is
f.g.  free and $H_*(\sigma^{-1}C)=0$.  In particular,
$L_{2i}(R,\sigma)$ is the Witt group of nonsingular $(-1)^i$-quadratic
$\sigma^{-1}R/R$-valued linking forms on f.g.  $\sigma$-torsion
$R$-modules of type ${\rm coker}(s:R^k \to R^k)$ ($s \in \sigma$),
and $\partial:L_{2i}(\sigma^{-1}R) \to L_{2i}(R,\sigma)$ is given by the
boundary construction for $\sigma^{-1}R$-nonsingular $(-1)^i$-quadratic
forms over $R$.

Given a ring $A$ let $\Pi^{-1}A[z,z^{-1}]$ be the universal
localization of $A[z,z^{-1}]$ inverting the set $\Pi$ of all
$A$-invertible square matrices over $A[z,z^{-1}]$. The canonical
ring morphism $A[z,z^{-1}] \to \Pi^{-1}A[z,z^{-1}]$ is an
injection with the universal property that every morphism of rings
$A[z,z^{-1}] \to R$ sending matrices in $\Pi$ to invertible
matrices over $R$ has a unique factorization $A[z,z^{-1}] \to
\Pi^{-1}A[z,z^{-1}] \to R$.

\begin{exa} {\rm For commutative $A$ $\Pi^{-1}A[z,z^{-1}]=P^{-1}A[z,z^{-1}]$
is the commutative localization of $A[z,z^{-1}]$ inverting the set
$P$ of all polynomials $p(z) \in A[z,z^{-1}]$ with $p(1) \in A$ a
unit.}\hfill$\qed$
\end{exa}

An involution on the ring $A$ is extended to the rings
$A[z,z^{-1}]$, $\Pi^{-1}A[z,z^{-1}]$ by $\overline{z}=z^{-1}$, as
before. As in Proposition 32.6 of Ranicki \cite{RHK}, the dual of
a Blanchfield $A[z,z^{-1}]$-module $B$ is given up to natural
isomorphism by
$$B\widehat{~}~=~{\rm Hom}_{A[z,z^{-1}]}(B,\Pi^{-1}A[z,z^{-1}]/A[z,z^{-1}])~.$$
An $A[z,z^{-1}]$-module morphism $\phi:B \to B\widehat{~}$ is the
same as a pairing
$$\phi~:~B \times B \to  \Pi^{-1}A[z,z^{-1}]/A[z,z^{-1}]$$
such that for all $x,x',y,y' \in B$, $a,b \in A[z,z^{-1}]$
$$\begin{array}{l}
\phi(x+x',y)~=~\phi(x,y)+\phi(x',y)~,\\[1ex]
\phi(x,y+y')~=~\phi(x,y)+\phi(x,y')~,\\[1ex]
\phi(ax,by)~=~b \phi(x,y) \overline{a} \in
\Pi^{-1}A[z,z^{-1}]/A[z,z^{-1}]~.
\end{array}$$

The quadratic $L$-groups $L_*$ of the Laurent polynomial extension
$A[z,z^{-1}]$ of a ring with involution $A$ split as
    $$L_{n+1}(A[z,z^{-1}])~=~L_{n+1}(A) \oplus L^p_n(A)$$
with $L^p_*$ the projective quadratic $L$-groups (Novikov
\cite{Nov}, Ranicki \cite{RL}).  The relative $L$-group
$L_{2i+2}(A[z,z^{-1}],\Pi)$ in the  localization exact sequence
    $$\begin{array}{ll}
   \dots \to L_{2i+2}(A[z,z^{-1}])& \to
    L_{2i+2}(\Pi^{-1}A[z,z^{-1}])\\[1ex]
    &\to L_{2i+2}(A[z,z^{-1}],\Pi) \to L_{2i+1}(A[z,z^{-1}]) \to \dots
    \end{array}$$
is the Witt group of nonsingular
$(-1)^{i+1}$-symmetric Blanchfield forms $(B,\phi)$ over $A[z,z^{-1}]$
such that $B$ admits a 1-dimensional f.g. free $A[z,z^{-1}]$-module resolution
$$0 \to C_1 \to C_0 \to B \to 0~.$$
As in Chapter 31 of \cite{RHK} let ${\rm LIso}^{2i}_p(A)$ be the Witt
group of nonsingular $(-1)^i$-symmetric Seifert forms over $A$.

\begin{thm} \label{cover3} {\rm (\cite[Prop 32.11]{RHK})}
The covering construction (\ref{Bcover}) defines an isomorphism
    $$B~:~{\rm LIso}^{2i}_p(A)\to L_{2i+2}(A[z,z^{-1}],\Pi)~;~
    (P,e,\theta) \mapsto B(P,e,\theta)~,$$
with Theorem \ref{cover2} giving an explicit inverse $B^{-1}$.
\hfill$\qed$
\end{thm}

The isomorphism $B^{-1}$ of \ref{cover3} is a generalization
of the projection
$$B~:~K_1(A[z,z^{-1}]) \to K_0(A)$$
of Bass, Heller and Swan \cite{BHS} and the projection
$$B~:~L_{2i+1}(A[z,z^{-1}]) \to L^p_{2i}(A)$$
of \cite{Nov} and \cite{RL} (where $B$ denotes Bass rather than
Blanchfield).

\begin{exa} {\rm The high-dimensional knot cobordism groups
$k:S^{2i-1} \subset S^{2i+1}$ ($i \geqslant 2$) are
    $$C_{2i-1}~=~{\rm LIso}^{2i}(\Z)~=~L_{2i+2}(\Z[z,z^{-1}],P)~.$$
See Chapters 33, 40 and 41 of \cite{RHK} for a more detailed
discussion. \hfill $\qed$}
\end{exa}

\begin{rem} \label{decomp5}
{\rm Theorems \ref{cover2}, \ref{cover3} give a new proof of the
result that every nonsingular $\eta$-symmetric Blanchfield form
$(B,\phi)$ over $A[z,z^{-1}]$ is isomorphic to the covering
$B(P,e,\theta)$ of a nonsingular $(-\eta)$-symmetric Seifert form
$(P,e,\theta)$ over $A$, with a corresponding isomorphism in the
Witt groups. For $A=\Z$, $\eta = \pm 1$ this was proved by a
variety of geometric and algebraic methods by Kearton \cite{K},
Levine \cite{L}, Trotter \cite{T} and Farber \cite{F}. For
arbitrary $A$ this was proved in Proposition 32.10 of Ranicki
\cite{RHK} using algebraic transversality for quadratic Poincar\'e
complexes over $A[z,z^{-1}]$. The novelty  is the explicit
algorithm for constructing $(P,e,\theta)$ from
$(B,\phi)$.}\vskip-4mm
$$\eqno{\qed}$$
\end{rem}

The expression of the Witt groups of $(-1)^{i+1}$-symmetric
Blanchfield forms over $A[z,z^{-1}]$ as the relative $L$-group
$L_{2i+2}(A[z,z^{-1}],\Pi)$ in the exact sequence of $L$-groups
 $$\begin{array}{ll}
   \dots \to L_{2i+2}(A[z,z^{-1}])& \to
    L_{2i+2}(\Pi^{-1}A[z,z^{-1}])\\[1ex]
    &\to L_{2i+2}(A[z,z^{-1}],\Pi) \to
    L_{2i+1}(A[z,z^{-1}]) \to \dots
    \end{array}$$
can be refined to an even more useful expression by inverting $1-z
\in A[z,z^{-1}]$. Write
$$A_z~=~A[z,z^{-1}]~,~A_{z,1-z}~=~A[z,z^{-1},(1-z)^{-1}]~.$$
For an $A$-module $P$ and an $A_z$-module $Q$ write
$$P_z~=~A_z\otimes_AP~,~P_{z,1-z}~=~A_{z,1-z}\otimes_AP~,~
Q_{1-z}~=~A_{z,1-z}\otimes_{A_z}Q~.$$
The element
$$s~=~(1-z)^{-1} \in A_{z,1-z}$$
is such that $s+\overline{s}=1 \in A_{z,1-z}$, so there is no
difference between $\pm$-quadratic and $\pm$-symmetric structures
(= forms, algebraic Poincar\'e complexes, $L$-groups) over
$A_{z,1-z}$. The cartesian square of rings with involution
$$\xymatrix{A_z \ar[r] \ar[d] & A_{z,1-z} \ar[d]\\
\Pi^{-1}A_z \ar[r] & \Pi^{-1}A_{z,1-z}}$$ induces excision
isomorphisms of relative $L$-groups
$$L_*(A_z,\Pi)~ \cong~ L_*(A_{z,1-z},\Pi)$$
and there is a commutative braid of exact sequences \vskip5mm

    $$\xymatrix@!C@C-70pt@R-10pt{
    L_{2i+3}(A_z,\Pi) \ar[dr]\ar@/^2pc/[rr] &&L_{2i+2}(A_{z,1-z})  \ar[dr]
    \ar@/^2pc/[rr]&& L_{2i+2}(A_z,(1-z)^{\infty})
    \ar[dr]\ar@/^2pc/[rr]&& L_{2i+1}(\Pi^{-1}A_z)\\
    & L_{2i+2}(A_z) \ar[ur]\ar[dr] &&
    L_{2i+2}(\Pi^{-1}A_{z,1-z})      \ar[ur]\ar[dr]    &&
    L_{2i+1}(A_z)\ar[ur]\ar[dr] & \\
    L_{2i+3}(A_z,(1-z)^{\infty})\ar[ur]\ar@/_2pc/[rr] &&
    L_{2i+2}(\Pi^{-1}A_z)\ar[ur]\ar@/_2pc/[rr] && L_{2i+2}(
    A_z,\Pi)\ar[ur]\ar@/_2pc/[rr] && L_{2i+1}(A_{z,1-z})}$$
\vskip10mm

\noindent
See Chapter 36 of \cite{RHK} for the identification of $L_{2i+2}(A_{z,1-z})$
with the Witt group of almost $(-1)^{i+1}$-symmetric forms $(P,\phi)$ over $A$,
with $P$ a f.g. free $A$-module and $\phi:P \to P^*$ an isomorphism such that
$1+(-1)^i(\phi^*)^{-1}\phi:P \to P$ is nilpotent (cf. Clauwens \cite{CL}).

\begin{thm} \label{invert}
The map $L_{2i+2}(A_z,\Pi) \to L_{2i+1}(A_{z,1-z})$ is 0, so that
    $$L_{2i+2}(A_z,\Pi)~=~{\rm coker}(L_{2i+2}(A_{z,1-z}) \to
    L_{2i+2}(\Pi^{-1}A_{z,1-z}))$$
The Witt class of the covering $B(P,e,\theta)$ of a nonsingular
$(-1)^i$-symmetric Seifert form $(P,e,\theta)$ over $A$ is
the Witt class of the nonsingular $(-1)^{i+1}$-quadratic form
$(P_{z,1-z},(1-z)\theta)$ over
$A_{z,1-z}$, modulo the indeterminacy coming from the
$(-1)^{i+1}$-quadratic Witt group of $A_{z,1-z}$.
\end{thm}
\begin{proof} Let $R$ be a ring with involution.
A 1-dimensional $(-1)^i$-quadratic Poincar\'e complex $(C,\psi)$
over $R$ with $\psi_0:C^0 \to C_1$ an isomorphism
(and $\widetilde{\psi}_0=0:C^1 \to C_0$) is null-cobordant
$$(C,\psi)~=~0 \in L_1(R,(-1)^i)~=~L_{2i+1}(R)~,$$
with a null-cobordism $(f:C \to D,(\delta\psi,\psi))$ defined by
    $$f~=~1~:~C_1 \to D_1~=~C_1~,~D_i~=~0~(i \neq 1)~,~\delta
    \psi~=~0~.$$
The nonsingular $(-1)^i$-quadratic formation corresponding to
$(C,\psi)$ is the boundary $\partial (C^0,\psi_1)$ of the
$(-1)^{i+1}$-quadratic form  $(C^0,\psi_1)$ over $R$.\\
\indent Now suppose that $R \to \sigma^{-1}R=S$ is an injective noncommutative
localization of rings with involution, so that there is defined a
localization exact sequence
    $$\dots \to L_{2i+2}(R) \to L_{2i+2}(S)
    \xymatrix@C-7.5pt{\ar[r]^{\partial}&} L_{2i+2}(R,\sigma)
    \to L_{2i+1}(R) \to \dots$$
with $L_{2i+2}(R,\sigma)$ the cobordism group of f.g. free 1-dimensional
$(-1)^i$-quadratic Poincar\'e complexes $(C,\psi)$ over $R$ such
that $1\otimes d:S\otimes_RC_1 \to S\otimes_RC_0$ is an $S$-module
isomorphism. If $(C,\psi)$ is such that $\psi_0:C^0 \to C_1$ is an
$R$-module isomorphism then (as above) $(C,\psi)=0 \in L_{2i+1}(R)$, and
$$1\otimes (\psi_1+(-1)^{i+1}\psi_1^*)~=~-1\otimes d\psi_0~:~
    S\otimes_RC^0 \to S\otimes_RC_0$$
is an $S$-module isomorphism, so that $S\otimes_R(C^0,\psi_1)$ is a nonsingular
$(-1)^{i+1}$-quadratic form over $S$ such that
    $$\begin{array}{l}
    (C,\psi)~=~ S\otimes_R(C^0,\psi_1)\\[1ex]
    \hskip25pt
    \in {\rm ker}(L_{2i+2}(R,\sigma) \to L_{2i+1}(R))~=~{\rm
    coker}(L_{2i+2}(R) \to L_{2i+2}(S))~.
    \end{array}$$
 \indent In particular,
if $(B,\phi)$ is a nonsingular $(-1)^{i+1}$-symmetric Blanchfield
form over $A_z$ then by Theorem \ref{cover2}
$(B,\phi)=B(P,e,\theta)$ is the covering of a nonsingular
$(-1)^i$-symmetric Seifert form $(P,e,\theta)$ over $A$. The
1-dimensional $(-1)^i$-quadratic Poincar\'e complex $(C,\psi)$
over $A_z$ defined by
$$\begin{array}{l}
d~=~\theta+(-1)^iz^{-1}\theta^*~:~C_1~=~P_z \to
C_0~=~P^*_z~,\\[1ex]
\psi_0~=~ 1-z~:~C^0~=~P_z \to C_1~=~P_z~,\\[1ex]
\psi_1~=~-(1-z)\theta~:~C^0~=~P_z \to C_0~=~P^*_z
\end{array}$$
has $1\otimes \psi_0:(C^0)_{1-z} \to (C_1)_{1-z}$  an
$A_{z,1-z}$-module isomorphism, so that
$$(B,\phi)_{1-z}~=~(C,\psi)_{1-z}~=~0 \in
    L_1(A_{z,1-z},(-1)^i)~=~L_{2i+1}(A_{z,1-z})~.$$
The nonsingular $(-1)^i$-quadratic formation over $A_{z,1-z}$
corresponding to $(C,\psi)_{1-z}$ is the boundary of the
$\Pi^{-1}A_{z,1-z}$-nonsingular $(-1)^{i+1}$-quadratic form
$$(C^0,\psi_1)_{1-z}~=~(P_{z,1-z},(1-z)\theta)$$
and
$$\begin{array}{ll}
(B,\phi)~=~(C,\psi)&=~(C,\psi)_{1-z}~=~
\partial (P_{z,1-z},(1-z)\theta)\\[1ex]
&\in {\rm ker}(L_{2i+2}(A_z,\Pi) \to L_{2i+1}(A_{z,1-z}))\\[1ex]
&\hskip10pt
=~{\rm ker}(L_{2i+2}(A_{z,1-z},\Pi) \to L_{2i+1}(A_{z,1-z}))\\[1ex]
&\hskip10pt
=~{\rm coker}(L_{2i+2}(A_{z,1-z}) \to L_{2i+2}(\Pi^{-1}A_{z,1-z}))~.
\end{array}$$
\vskip-3mm
\end{proof}

\begin{exa} {\rm The expression for the Witt group of Blanchfield
forms given by Theorem \ref{invert} in the case $A=\Z$ gives the
following expression for the cobordism class of a high-dimensional knot.
Let $k:S^{2i-1} \subset S^{2i+1}$ be a knot with exterior
$$(M^{2i+1},\partial M)~=~({\rm cl.}(S^{2i+1} \backslash
 k(S^{2i-1}) \times D^2),S^{2i-1} \times S^1)$$
and Seifert surface $N^{2i} \subset S^{2i+1}$. Keeping
$\partial N=k(S^{2i-1})$ fixed push $N$ into the interior
of $D^{2i+2}$ to obtain a codimension 2 embedding $N \subset D^{2i+2}$
with trivial normal 2-plane bundle. The exterior is a $(2i+2)$-dimensional
manifold with boundary
    $$(L^{2i+2},\partial L)~=~({\rm cl.}(D^{2i+2} \backslash N \times D^2),
    M \cup_{\partial M} N \times S^1)~.$$
Assume that $\pi_j(M) \cong \pi_j(S^1)$ for $1 \leqslant j \leqslant
i-1$ (as may be arranged by surgery below the middle dimension),
so that $N$ can be chosen to be $(i-1)$-connected, and $\pi_j(L)
\cong \pi_j(S^1)$ for $1 \leqslant j \leqslant i$.  As in Proposition
27.8 of \cite{RHK} there is defined an $(i+1)$-connected
$(2i+2)$-dimensional normal map of triads
$$\begin{array}{l}
(f,b)~:~(L;M,N \times S^1;S^{2i-1} \times S^1)\\[1ex]
\hskip50pt \to (D^{2i+2} \times  [0,1];D^{2i+2} \times \{0\},
D^{2i+2} \times \{1\};k(S^{2i-1})\times [0,1]) \times S^1
\end{array}$$
with target a $(2i+2)$-dimensional geometric Poincar\'e triad. The
nonsingular $(-1)^i$-symmetric Seifert form $(H_i(N),e,\theta)$
over $\Z$ determines the kernel $(-1)^{i+1}$-quadratic form over
$\Z_z=\Z[z,z^{-1}]$
$$(K_{i+1}(\overline{L}),\psi)~=~(H_i(N)_z,(1-z)\theta)$$
(cf.  Ko \cite{Ko}, Cochran, Orr and Teichner \cite{COT}).  For $i
\geqslant 2$ the knot cobordism class of $k$ is the Witt class of
$(H_i(N),e,\theta)$, or equivalently the Witt class of the
nonsingular $(-1)^{i+1}$-symmetric Blanchfield form
$(H_i(\overline{M}),\phi)=B(H_i(N),e,\theta)$ over $\Z_z$. Theorem
\ref{invert} identifies the knot cobordism class with the Witt
class (modulo the indeterminacy) of the induced nonsingular
$(-1)^{i+1}$-quadratic form over
$\Z_{z,1-z}=\Z[z,z^{-1},(1-z)^{-1}]$
$$\begin{array}{ll}
[k]=(H_i(N),e,\theta)=(H_i(\overline{M}),\phi)=
(K_{i+1}(\overline{L}),\psi)_{1-z}=(H_i(N)_{z,1-z},(1-z)\theta)\\[1ex]
\in C_{2i-1}={\rm LIso}^{2i}(\Z)=L_{2i+2}(\Z_z,P)= {\rm coker}(
L_{2i+2}(\Z_{z,1-z}) \to L_{2i+2}(P^{-1}\Z_{z,1-z}))
\end{array}$$
with $P=\{p(z)\vert p(1)=1\} \subset \Z_z$ the multiplicative
subset of Alexander polynomials. See Proposition 36.3 of
\cite{RHK} for the computation of the indeterminacy
    $$L_{2i+2}(\Z_{z,1-z})~=~
    L^{2i+2}(\Z)~=~ \begin{cases} 0&\hbox{if $i \equiv 0(\bmod\,2)$}\\
    \Z~({\rm signature})&\hbox{if $i \equiv 1(\bmod\,2)$}~.
    \end{cases}$$
\hfill$\Box$}
\end{exa}

\providecommand{\bysame}{\leavevmode\hbox to3em{\hrulefill}\thinspace}

\end{document}